\documentclass{article}
\usepackage{latexsym}
\usepackage{amsmath}
\usepackage{amssymb}
\usepackage[dvips]{graphicx}

\newcommand{\en}{_{\varepsilon,n}}
\newcommand{\rtre}{{\mathbb R}^3}
\newcommand{\rdue}{{\mathbb R}^2}
\newcommand{\erre}{{\mathbb R}}

\newcommand{\graph}[1]{{\rm graph}\, #1}

\newcommand{\dive}{{\rm div}}
\newcommand{\I}[3]{I( \xi, #1 , #2, #3 )}
\newcommand{\emme}{\rho_{\varepsilon, n}}
\newcommand{\emmee}{\rho_{\varepsilon}}
\newcommand{\thetae}{\theta_{\varepsilon}}
\newcommand{\ren}{r_{\varepsilon}^h}

\newcommand{\thetak}{\theta_{\varepsilon, n_{\varepsilon}}}
\newcommand{\thetat}{\tilde{\theta}}
\newcommand{\stheta}{\overline{\theta}_{\varepsilon}^h}
\newcommand{\qed}{\hfill $\Box$}
\newcommand{\diff}{\partial}
\newcommand{\varx}{\tilde{x}}
\newcommand{\vary}{\tilde{y}}
\newcommand{\ddxi}{\tau_{\xi}}
\newcommand{\ddeta}{\tau_{\eta}}
\newcommand{\ddz}{e_{z}}
\newcommand{\haus}{{\cal H}^{1}} 
\newcommand{\Om}{\Omega}
\newcommand{\curv}{{\rm curv}\,}
\newcommand{\nablav}{\nabla}
\newcommand{\usec}{\overline{U'}}
\newcommand{\bclaim}[2]{\vspace{.3cm}\noindent {\bf Claim #1.} {\sl #2}}
\newcommand{\eclaim}[1]{\vspace{.1cm}\noindent {\sc Proof of the claim.} #1
\vspace{.1cm}}  
\newcommand{\bstep}[2]{\vspace{.5cm}\noindent {\bf STEP #1.} {  #2}\vspace{.5cm}}
\newcommand{\vuu}{v_1}
\newcommand{\vud}{v_2}
\newcommand{\vui}{v_i}
\newcommand{\vuup}{\tilde{v}_1}
\newcommand{\vudp}{\tilde{v}_2}
\newcommand{\vuip}{\tilde{v}_i}
\newcommand{\vuut}{\hat{v}_1}
\newcommand{\vudt}{\hat{v}_2}
\newcommand{\vuit}{\hat{v}_i}
\newcommand{\aplimsup}{{\rm ap}\limsup}

\makeatletter
\@addtoreset{equation}{section}
\makeatother
\newtheorem{theorem}{Theorem}[section]
\newtheorem{lemma}[theorem]{Lemma}

\newtheorem{prop}[theorem]{Proposition}
\newtheorem{cor}[theorem]{Corollary}
\newtheorem{definition}[theorem]{Definition}

\begin{document}

\null
\vspace{1cm}
\begin{center}
{\LARGE Local calibrations for minimizers of the  Mumford-Shah 

\vspace{.5cm}

functional with a regular discontinuity set}
\vspace{1cm}

{\normalsize Maria Giovanna \sc Mora}

\vspace{.25cm}
{\normalsize Massimiliano \sc Morini}

\vspace{.5cm}

{\normalsize S.I.S.S.A.}\\
{\normalsize via Beirut 2-4, 34014 Trieste, Italy}\\
{\normalsize e-mail: \tt mora@sissa.it, 
morini@sissa.it}

\vspace{2cm}

\begin{minipage}[t]{11.5cm}
{\footnotesize
\centerline{ {\bf Abstract}}
\noindent Using a calibration method, 
we prove that, if $w$ is a function which satisfies
all Euler conditions for the Mumford-Shah functional on a two-dimensional
open set $\Omega$, and the discontinuity set $S_w$ of $w$ is a regular curve 
connecting two boundary points, then there exists a uniform
neighbourhood $U$ of $S_w$ such that
$w$ is a minimizer of the Mumford-Shah functional on $U$ with respect 
to its own boundary conditions on $\partial U$.
We show that Euler conditions do not guarantee in general the minimality
of $w$ in the class of functions with the same boundary value of $w$ 
on $\partial \Omega$ and whose
extended graph is contained in a neighbourhood of the extended graph of $w$,
and we give a sufficient condition in terms of the geometrical properties
of $\Omega$ and $S_w$ under which this kind of minimality holds.}
\vspace{.5cm}\\


\vspace{.5cm}
\end{minipage}
\end{center}

\setcounter{page}{-1}
\thispagestyle{empty}

\vfill

\begin{center}
{\normalsize Ref. S.I.S.S.A. 21/2000/M  (March 2000)}
\end{center}

\pagebreak

\clearpage
\phantom{a}\thispagestyle{empty}
\pagebreak

\section{Introduction}

This paper deals with local minimizers of the Mumford-Shah functional
(see \cite{Mum-Sha1} and \cite{Mum-Sha2})
\begin{equation} \label{g1}
\int_{\Om} |\nabla u(x,y)|^{2} dx\,dy +  \haus(S_{u}) \,,
\end{equation}
where $\Omega$ is a bounded open subset of $\rdue$ with a Lipschitz 
boundary, $\haus$ is the one-dimensional Hausdorff measure,
$u$ is the unknown function in the space 
$SBV(\Om)$ of special functions of bounded 
variation in $\Om$,
$S_{u}$ is the set of essential 
discontinuity points of $u$, while $\nabla u$ 
denotes its approximate gradient (see \cite{Amb} or \cite{Amb-Fus-Pal}).
\begin{definition}
We say (as in \cite{Alb-Bou-DM}) that $u$ is a Dirichlet minimizer of (\ref{g1}) in $\Om$ 
if it belongs to $SBV(\Om)$ and satisfies the inequality
$$\int_{\Om} |\nabla u(x,y)|^{2} dx \,dy + \haus(S_{u}) \,\le\,
\int_{\Om} |\nabla v(x,y)|^{2} dx \,dy + \haus(S_{v})$$
for every function $v\in SBV(\Om)$ with the same trace as $u$ on $\partial\Om$. 
\end{definition}
Suppose that $u$ is a Dirichlet minimizer 
of (\ref{g1}) in $\Om$ and that $S_{u}$ is a regular curve. Then the following 
equilibrium conditions are satisfied (see \cite{Mum-Sha1} and
\cite{Mum-Sha2}):
\begin{description}
\item[i)] $u$ is harmonic on $\Om\setminus S_{u}$;
\item[ii)] the normal derivative of $u$ vanishes on both sides of $S_{u}$;
\item[iii)] the curvature of $S_{u}$ is equal to the
difference of the squares of the tangential derivatives of $u$
on both sides of $S_{u}$.
\end{description}
Elementary examples show that conditions i), ii), and iii) are not 
sufficient for the Dirichlet minimality of $u$.

In this paper we prove that, if 
$S_u$ is an analytic curve connecting two points of $\partial \Om$,
then i), ii), iii) are also sufficient for the Dirichlet minimality
of $u$ in small domains.
In other words, for every $(x_0, y_0)$ in $\Om$, there is
an open neighbourhood $U$ of $(x_0, y_0)$ such that
$u$ is a Dirichlet minimizer of (\ref{g1}) in $U$. If $(x_0,y_0)$ does not lie
on $S_u$, this fact is well known and can be proved by the calibration method
(see \cite{Alb-Bou-DM}); so  
the interesting case is when we consider points belonging to $S_u$:  in this situation 
we have a stronger result, since we can prove that the Dirichlet minimality actually holds 
in a uniform neighbourhood of the discontinuity set.
The analyticity assumption for $S_u$ does not seem too restrictive: it has
been proved that the regular part of the discontinuity set of a minimizer is
of class $C^{\infty}$ and 
it is a conjecture that it is analytic
(see \cite{Amb-Fus-Pal}).

Let us give the precise statement of the result.
\begin{theorem}\label{teoparz}
Let $\Om_0$ be a connected open subset of $\rdue$ and
$\Gamma$ be a simple analytic curve in $\Om_0$ connecting two points 
of the boundary. Let $u$ be a function in $H^1( \Om_0 \setminus \Gamma)$ with
$S_u = \Gamma$, with different traces at every point of $\Gamma$, and
satisfying the Euler conditions i), ii), and iii) in $\Om_0$ (for the precise
formulation of these conditions, see Section 2). Finally, let $\Om$ be an open
set with Lipschitz boundary,  compactly contained in $\Om_0$,
such that $\Om \cap \Gamma \neq \emptyset$.
Then there exists an open neighbourhood $U$ of $\Gamma \cap \overline{\Om}$
contained in $\Om_0$ such that $u$ is a Dirichlet minimizer in $U$ of
the Mumford-Shah functional (\ref{g1}).
\end{theorem}
This theorem generalizes the result of Theorem 4.2
of \cite{DM-Mor-Mor} in two directions: the discontinuity set $S_u$
can be any analytic curve and the Dirichlet minimality of $u$ is proved
in a uniform neighbourhood of $S_u\cap \overline{\Om}$.
The proof is obtained, as in \cite{DM-Mor-Mor}, by the calibration
method introduced in \cite{Alb-Bou-DM}. The original idea 
of the new construction essentially
relies on the definition of the calibration around the graph of $u$:
here it is obtained using the gradient field
of a family of harmonic functions, whose graphs fiber a
neighbourhood of the graph of $u$. This technique seems to have some
similarities with the classical method of the Weierstrass fields, where
the proof of the minimality of a candidate $u$ is obtained by 
the construction of
a slope field starting from a family of solutions of 
the Euler equation, whose graphs foliate a neighbourhood of the graph
of $u$.

In this paper
we are also interested in a different type of minimality:
in Theorem \ref{teoparz} we compare $u$ with perturbations which can be very large, 
but concentrated in a fixed small domain; we wonder 
if a minimality property is preserved also
when we admit as competitors perturbations of $u$ with $L^{\infty}$-norm very small 
outside a small neighbourhood of $S_u$,
but support possibly coinciding with $\overline{\Omega}$.

This is made precise by the following definition.

\begin{definition} \label{graphmin}
A function $u \in SBV(\Omega)$ is a local 
graph-minimizer 
in $\Omega$ if there exists a suitable neighbourhood $U$ of the extended graph $\Gamma_u$ 
of $u$ (for the notion of extended graph, see Section 2) such that
$$\int_{\Om} |\nabla u(x,y)|^2 dx \, dy + {\cal H}^{1}(S_u)
\leq \int_{\Om} |\nabla v(x,y)|^2 dx \, dy + {\cal H}^{1}(S_v)$$ 
for every $v\in SBV(\Omega)$ with the same trace as $u$ on $\partial \Omega$ 
and whose extended graph $\Gamma_v$ is contained in $U$.
\end{definition}
In \cite{Alb-Bou-DM} it is proved that any harmonic function 
defined on $\Omega$ is a local graph-minimizer
whatever $\Omega$ is. If the function presents some discontinuities, 
what we discover is that the graph-minimality may fail when $\Omega$ is too large, 
even in the case of rectilinear discontinuities, as the counterexample 
given in Section 4 shows.

To get the graph-minimality we have to add some restrictions on the domain $\Omega$. 
To this aim we introduce a suitable quantity which seems useful to describe
the right geometrical interaction between $S_u$ and $\Omega$. Given an open
set $A$ (with Lipschitz boundary) and a portion $\Gamma$ of
$\partial A$ (with nonempty relative interior in $\partial A$), we define
$K(\Gamma,A)$  by the variational problem
\begin{equation}\label{capacita}
K(\Gamma,A): =
\inf \left\{ \int_{A} |\nabla v(x,y)|^2 dx \, dy :
\ v\in H^1(A), \ \int_{\Gamma} v^2 d{\cal H}^1=1, 
\hbox{ and }v=0 \hbox{ on }\partial A\setminus\Gamma \right\}.
\end{equation}
First of all, it is easy to see that 
in the problem above the infimum is attained; moreover, 
the notation is well chosen since $K(\Gamma,A)$ is a quantity  
depending only on $\Gamma$ and $A$, which describes 
a kind of ``capacity'' of the prescribed portion of the boundary 
with respect to the whole open set. 
Note also that if $A_1\subset A_2$, and $\Gamma_1 \subset \Gamma_2$,
then $K(\Gamma_1, A_1)\geq K(\Gamma_2, A_2)$, 
which suggests that if $K(\Gamma,A)$ is very large, 
then $A$ is thin in some sense.
It is convenient to give the following definition.

\begin{definition}\label{gadm}
Given a simple analytic curve $\Gamma$, we say that
an open set $\Omega$ is $\Gamma$-admissible if
it is bounded, $\Gamma\cap \overline{\Omega}$ connects two points of
$\partial \Omega$, and $\Omega\setminus\Gamma$ has two connected components,
which have Lipschitz boundary. 
\end{definition}

The following theorem gives a sufficient condition
for the graph-minimality in terms of $K(\Gamma,\Omega)$  and of the geometrical
properties of the curve. We denote the length of $\Gamma$ by
$l(\Gamma)$, its curvature by $\curv \Gamma$, and the
$L^{\infty}$-norm of $\curv\Gamma$ by $k(\Gamma)$.

\begin{theorem}\label{teoglob}
Let $\Omega_0$, $\Omega$, $u$, and $\Gamma=S_u$ satisfy 
the same assumptions as in Theorem \ref{teoparz}; 
suppose that $\Omega$ is $\Gamma$-admissible and denote by
$\Omega_1$ and $\Omega_2$ the two connected components of
$\Omega\setminus\Gamma$, by $u_i$ the restriction of $u$ to $\Omega_i$, and by
$\partial_{\tau} u_i$ its tangential derivative on $\Gamma$. There exists an
absolute constant $c>0$ (independent of $\Omega_0$, $\Omega$, $\Gamma$, and
$u$) such that if 
\begin{equation}\label{sufficiente}
 \frac{\min_{i=1,2}K(\Gamma\cap\Omega,\Omega_i)}
{1 + l^2(\Gamma\cap \Omega) + l^2(\Gamma\cap \Omega) k^2(\Gamma\cap
\Omega)} > c  \sum_{i=1}^2 \|\partial_{\tau}
u_i\|^2_{C^1(\Gamma\cap\Omega)},
\end{equation}
then $u$ is a local graph-minimizer on $\Omega$.
\end{theorem}
Remark that condition (\ref{sufficiente}) imposes a restriction 
on the size of $\Omega$ depending on the behaviour of $u$ along $S_u$: 
if $u$ has large or very oscillating tangential derivatives, we have to take
$\Omega$ quite small to guarantee that (\ref{sufficiente}) is satisfied. In
the special case of a locally constant function $u$, condition
(\ref{sufficiente}) is always fulfilled whatever the domain is; so $u$ is a
local graph-minimizer whatever $\Omega$ is, in agreement with a result that
will be proved in the final version of \cite{Alb-Bou-DM}. 

The plan of the paper is the following: in Section 2, we fix some notation
and recall the main result of \cite{Alb-Bou-DM}; Section 3 contains the
proof of Theorem \ref{teoparz}; finally, Section 4 is devoted to
the graph-minimality: we give a counterexample when (\ref{sufficiente}) is
violated, we prove Theorem \ref{teoglob}, and present some qualitative
properties of $K(\Gamma,\Omega)$.  

\section{Preliminary results}

Given any subset $A$ of $\rdue$ and $\delta>0$,
we denote by $A_{\delta}$ the $\delta$-neighbourhood of $A$,
defined by
$$A_{\delta} := \{ (x_0,y_0)\in \rdue : \, \exists (x,y)\in A \ \hbox{such
that}\
|(x-x_0, y-y_0)|< \delta \}.$$
Let $\Om$ be an open set in $\rdue.$
If $v\in SBV(\Om)$, for every $(x_0,y_0)\in \Om$
we put
$$v^+(x_0,y_0):=  {\aplimsup}_{(x,y)\to (x_0,y_0)} v(x,y)
\qquad \hbox{and} \qquad
v^-(x_0,y_0):= {{\rm ap}\liminf}_{(x,y)\to (x_0,y_0)} v(x,y),$$
(see \cite{Amb-Fus-Pal}).
We recall that 
$v^+ = v^-$ ${\cal H}^{1}$-a.e.
in $\Om \setminus S_v$, while
for ${\cal H}^{1}$-a.e. $(x_0,y_0)\in S_v$
$$v^{\pm}(x_0,y_0) = \lim_{r\to 0^+} \frac{1}{{\cal L}^2(B^{\pm}_r(
x_0,y_0))}\int_{B^{\pm}_r(x_0,y_0)}
v(x,y)\,dx\,dy, $$
where $B^{\pm}_r(x_0,y_0)$ is the intersection of the ball of radius $r$
centred at $(x_0,y_0)$ 
with the half-space $\{ (x,y)\in \rdue: 
\pm(x-x_0, y-y_0)\cdot \nu_v(x_0,y_0) \geq 0 \}$,
where 
the vector $\nu_v(x_0,y_0)$ is the normal vector
to $S_v$ at $(x_0,y_0)$ (which is defined ${\cal H}^1$-a.e. on $S_v$).
The extended graph of $v$ is the set
$$\Gamma_v:= \{ (x,y,t)\in \Om { \times } \erre:
v^-(x,y)\leq t \leq v^+(x,y) \}.$$
Let $\Gamma$ be a smooth curve in $\Om$.
Fix an orientation of $\Gamma$ and call 
$\nu$ the corresponding normal vector field to $\Gamma$.
Let $\xi \mapsto
(x(\xi), y(\xi))$
be a parameterization of $\Gamma$ by the arc-length. The (signed)
curvature is defined by
\begin{equation}\label{curvatura1}
\curv\Gamma (\xi) = - (\ddot{x}(\xi), \ddot{y}(\xi))\cdot \nu(\xi);
\end{equation}
since the two vectors in (\ref{curvatura1}) are parallel, it follows that
\begin{equation}\label{curvq}
[\curv\Gamma (\xi)]^2 = (\ddot{x}(\xi))^2 + (\ddot{y}(\xi))^2.
\end{equation}

Let $u\in SBV(\Om)$ be a function such that 
$S_u = \Gamma$.
We say that $u$ satisfies the Euler conditions
for the Mumford-Shah functional in $\Om$ if
\begin{description}
	\item[i)]	$u$ is harmonic in $\Om \setminus \Gamma$ and $u\in
			H^1(\Omega \setminus \Gamma)$,

	\item[ii)]	$\displaystyle \frac{\diff u}{\diff \nu}=0$ on
			$\Gamma$,

	\item[iii)]	$|\nabla u^+|^2 - |\nabla u^-|^2 = \curv\Gamma$ at
			every point of $\Gamma$,
\end{description}
where $\nabla u^{\pm}$ denote the traces of $\nabla u$ on $\Gamma$.

If $U$ is any open subset of $\rtre$,
we shall consider the collection ${\cal F}(U)$ of all
piecewise $C^1$ vector fields $\varphi: U
\to \rdue { \times } \erre$ with the following property:
there exists a finite family $(A_i)_{i\in I}$ of
pairwise disjoint open subsets of $U$
such that the family of their closures covers $U$,
$\partial A_i \cap U$ is a Lipschitz
surface without boundary for every $i\in I$,
and $\varphi |_{A_i} \in C^1(\overline{A_i}, \rdue { \times } \erre).$

For every vector field $\varphi : 
U \to \rdue { \times } \erre$ we define the maps
$\varphi^{x}, \; \varphi^{y}, \; \varphi^{z} : U \to 
\erre$ by 
$$\varphi (x,y,z) = (\varphi^{x}(x,y,z), \varphi^{y}(x,y,z), 
\varphi^{z}(x,y,z)).$$

Let $U$ be an open neighbourhood of $\Gamma_u$ such that
the intersection with every straight vertical line
is connected.
A {\it calibration \/} for $u$ in $U$ is a bounded vector field $\varphi 
\in {\cal F}(U)$ which is continuous on the graph of $u$ and
satisfies the following properties:
\begin{description}
\item[(a)] $\dive \varphi = 0$ in the sense of distributions in 
$U$;
\item[(b)] $(\varphi^{x}(x,y,z))^{2} + (\varphi^{y}(x,y,z))^{2} \leq 
4 \varphi^{z}(x,y,z)$ at every continuity point $(x,y,z)$ of $\varphi$;
\item[(c)]  $(\varphi^{x}, \varphi^{y})(x,y,u(x,y)) = 2 \nabla  u(x,y)$
and
$\varphi^{z}(x,y,u(x,y)) = |\nabla u(x,y)|^{2}$ for every $(x,y)\in 
\Om \setminus S_u$;
\item[(d)] $\displaystyle \left( \int_{s}^{t} \varphi^{x}(x,y,z)\, dz 
\right)^{2} + \left( \int_{s}^{t} \varphi^{y}(x,y,z)\, dz 
\right)^{2} \leq 1$ for every $(x,y)\in \Om$ and for every $s, 
t$ such that $(x,y,s), (x,y,t)\in U$;
\item[(e)] $\displaystyle \int_{u^{-}(x,y)}^{u^{+}(x,y)} (\varphi^{x},
\varphi^{y})(x,y,z)\, dz =
\nu_u (x,y)$ 
for every $(x,y)\in S_u$.
\end{description}

The following theorem is proved in \cite{Alb-Bou-DM}.

\begin{theorem}
If there exists a calibration $\varphi$ for $u$ in $\Om { \times } \erre$, 
then $u$ is a Dirichlet minimizer of the Mumford-Shah 
functional (\ref{g1}) in $\Om$.
\end{theorem}
What the authors actually prove (but it is not explicitly
remarked), is the following more general statement.

\begin{theorem}
Let $U$ be an open neighbourhood of $\Gamma_u$ such that
the intersection with every straight vertical line
is connected. 
If there exists a calibration $\varphi$ for $u$ in $U$, 
then 
$$\int_{\Om} |\nabla u(x,y)|^2dx\,dy + {\cal H}^{1}(S_u\cap \Om)
\leq \int_{\Om} |\nabla v(x,y)|^2dx\,dy + {\cal H}^{1}(S_v)$$
for every $v\in SBV(\Om)$ such that $v=u$ on $\partial \Om$
and $\Gamma_v \subset U$.
\end{theorem}

\section{Proof of Theorem \ref{teoparz}}

\begin{lemma}\label{alberti}
Let $U$ be an open subset of $\rdue$ and 
$I$, $J$ be two real intervals. Let $u: U {\times} J \to I$
be a function of class $C^1$ such that
\begin{itemize}
\item	$u(\cdot, \cdot \, ; s)$ is harmonic for every $s\in J$;

\item	there exists a $C^1$ function $t :U { \times } I \to J$
	such that $u(x,y ; t(x,y;z ))=z$. 
\end{itemize}
Then, if we define in $U { \times } I$ the vector field 
$$\phi (x,y,z):= (2 \nabla u(x,y; t(x,y;z )), |\nabla u(x,y; t(x,y;z ))|^2),$$
where $\nabla u(x,y; t(x,y;z ))$ denotes the gradient of $u$ with respect to the variables
$(x,y)$ computed at $(x,y; t(x,y;z ))$, 
$\phi$ is divergence free in $U{ \times } I$.
\end{lemma}

\noindent {\sc Proof of the lemma.} 
Let us compute the divergence of $\phi$:
\begin{eqnarray}\label{laplac}
\dive \phi(x,y,z) & = & 2 \triangle
u(x,y; t(x,y;z)) + 2 \diff_s \nabla u(x,y;t(x,y;z))\cdot
\nabla t(x,y;z) \nonumber \\
& & + 2 \diff_z t(x,y;z)\,\nabla u(x,y;t(x,y;z)) 
\cdot \diff_s \nabla  u(x,y;t(x,y;z)),
\end{eqnarray}
where $\triangle u(x,y; t(x,y;z))$ denotes the laplacian of $u$ with respect to
$(x,y)$ computed at $(x,y; t(x,y;z))$, and $\nabla t(x,y;z)$ denotes the
gradient of $t$ with respect to $(x,y)$.  
By differentiating the identity
verified by the function $t$ first with respect to $z$ and with respect to
$(x,y)$, we derive that $$\diff_s u(x,y;t(x,y;z))\, \diff_z t(x,y;z) =1,
\qquad \nabla u(x,y;t(x,y;z)) + \diff_s u(x,y;t(x,y;z))
\, \nabla t(x,y;z)=0.$$
Using these identities and substituting
in (\ref{laplac}), we finally obtain
$$\dive \phi(x,y,z) = 2 \triangle u(x,y;t(x,y;z))=0,$$
since by assumption $u$ is harmonic with respect to $(x,y)$.
\qed

\

\noindent {\sc Proof of Theorem \ref{teoparz}.} 
In the sequel, the intersection $\Gamma \cap
\overline{\Om}$ will be still denoted by $\Gamma$.
Let
$$\Gamma: \begin{cases}
x=x(s) \\
y=y(s) \\
\end{cases}$$
be a parameterization by the arc-length,
where $s$ varies in $[0, l(\Gamma)]$;
we choose as orientation 
the normal vector field $\nu(s)=(-\dot{y}(s), \dot{x}(s))$.

By Cauchy-Kowalevski Theorem (see \cite{Joh}) there exist an open neighbourhood 
$U$ of $\Gamma$
contained in $\Om_0$ and a harmonic function $\xi$ defined on $U$ 
such that
$$\xi(\Gamma (s))=s \qquad \hbox{and} \qquad \frac{\diff \xi}{\diff
\nu}(\Gamma(s))=0.$$ We can suppose that $U$ is simply connected.
Let $\eta: U \to \rdue$ be the harmonic conjugate of $\xi$
that vanishes on $\Gamma$, i.e., the function satisfying
$\diff_x \eta(x,y) = - \diff_y \xi(x,y)$, 
$\diff_y \eta(x,y) = \diff_x \xi(x,y)$, and
$\eta(\Gamma(s))=0$.

Taking $U$ smaller if needed, we can suppose that the map
$\Phi(x,y):= (\xi(x,y), \eta(x,y))$ is invertible on $U$. We call
$\Psi$ the inverse function $(\xi, \eta) \mapsto (\varx(\xi,\eta), \vary(\xi,\eta))$, 
which is defined 
in the 
open set $V:=\Phi (U)$.
Note that, if $U$ is small enough, then $(\varx(\xi,\eta),\vary (\xi, \eta))$ 
belongs to $\Gamma$ 
if and only if $\eta=0$. Moreover, 
\begin{equation}\label{psijacob}
D\Psi = \left( \begin{array}{cc}
\diff_{\xi} \varx & \diff_{\eta} \varx  \\
\diff_{\xi} \vary  & \diff_{\eta} \vary  
\end{array} \right) =
\frac{1}{|\nabla \xi|^2} \left(
\begin{array}{cc}
\diff_x \xi & \diff_x \eta  \\
\diff_y \xi  & \diff_y \eta
\end{array}
\right),
\end{equation}
where, in the last formula, all functions are computed at $(x,y)= \Psi(\xi, \eta)$,
and so
\begin{equation}\label{001}
\diff_{\xi} \varx = \diff_{\eta} \vary\qquad \hbox{and}\qquad 
\diff_{\eta} \varx = - \diff_{\xi} \vary.
\end{equation} 
In particular, $\varx$ 
and $\vary$ are harmonic.

On $U$ we will use the coordinate system $(\xi,\eta)$ given by $\Phi$.
By (\ref{psijacob}) the canonical 
basis of the tangent space to $U$ at a point $(x,y)$
is given by
\begin{equation}\label{diff}
\ddxi = \frac{\nabla \xi}{|\nabla \xi|^2}, \qquad \ddeta =\frac{\nabla \eta}{|\nabla \eta|^2}.
\end{equation}
For every $(\xi,\eta)\in V$,
let $G(\xi,\eta)$ be the matrix associated with the first fundamental form
of $U$ in the coordinate system $(\xi,\eta)$,
and let $g(\xi,\eta)$  be its determinant. 
By (\ref{psijacob}) and (\ref{diff}), 
\begin{equation}\label{gi}
g = ((\diff_{\xi} \varx)^2 + (\diff_{\xi} \vary )^2)^{2} = \frac{1}{|\nabla \xi 
(\Psi) |^4}.
\end{equation}
We set $\gamma(\xi,\eta) = \sqrt[4]{g(\xi,\eta)}$.

From now on we will assume that $V$ is symmetric
with respect to
$\{(\xi,\eta)\in \Phi(U): \eta=0\}$.

Note that we can write the function $u$ in this new coordinate system as
$$u(\xi,\eta)= 
\begin{cases}
u_1 (\xi,\eta) & \text{if $(\xi,\eta)\in V$, $\eta <0$,} \\
u_2 (\xi,\eta) & \text{if $(\xi,\eta)\in V$, $\eta >0$,}
\end{cases}$$
where we can suppose that $u_1$ and $u_2$ are defined in $V$ (indeed, $u_1$ is a priori
defined only on the set $\{(\xi,\eta)\in V:\eta<0\}$, but it can be extended to $V$
by reflection; an analogous argument applies to $u_2$), 
$0< u_1(\xi, 0) < u_2(\xi, 0)$ for every $(\xi,0)\in V$,
and 
\begin{description}
\item[i)] $\diff^2_{\xi \xi} u_i (\xi,\eta) + \diff^2_{\eta \eta} u_i(\xi,\eta) =0$ for $i=1,2$;

\item[ii)] $\diff_{\eta} u_1 (\xi,0) = \diff_{\eta} u_2 (\xi, 0) = 0$;

\item[iii)] $(\diff_{\xi} u_2 (\xi, 0))^2 - (\diff_{\xi} u_1 (\xi, 0))^2=
		\curv \Gamma(\xi)$.
\end{description}

The calibration $\varphi(x,y,z)$ on $U { \times } \erre$ 
will be written as 
\begin{equation}\label{calibra2}
\varphi(x,y,z) = \frac{1}{\gamma^2 (\xi(x,y), 
\eta(x,y))}\phi(\xi(x,y), 
\eta(x,y), z),
\end{equation}
where $\phi: V { \times } \erre \to \rtre$ can be represented by
\begin{equation}\label{calibra}
\phi (\xi,\eta,z)=
\phi^{\xi} (\xi,\eta,z) \ddxi +  \phi^{\eta} (\xi,\eta,z) \ddeta +  \phi^z (\xi,\eta,z) 
\ddz,
\end{equation}
where $\ddz$ is the third vector of the canonical basis of $\rtre$, 
and $\ddxi$, $\ddeta$ are computed at the point $\Psi (\xi,\eta)$.
We now reformulate the conditions of Section 2 
in this new coordinate system.
It is known from Differential Geometry (see, e.g., \cite[Proposition 
3.5]{Cha}) that,
if $X = X^{\xi} \ddxi + X^{\eta} \ddeta$ 
is a vector field on $U$,
then the divergence of $X$ is given by
\begin{equation}\label{div}
\dive X = \frac{1}{\gamma^2} (\diff_{\xi}(\gamma^2 X^{\xi}) +
\diff_{\eta} (\gamma^2 X^{\eta})).
\end{equation}
Using (\ref{diff}), (\ref{gi}), (\ref{calibra2}), (\ref{calibra}), and
(\ref{div})
it turns out
that $\varphi$ is a calibration if
the following conditions are satisfied: 
\begin{description}\label{calk}
\item[(a)] $\displaystyle 
     \diff_{\xi} \phi^{\xi} +
	\diff_{\eta} \phi^{\eta} +
	\diff_{z} \phi^z = 0$ for every $(\xi,\eta,z) \in V { \times } \erre$;
\item[(b)] $\displaystyle (\phi^{\xi} (\xi,\eta,z))^2 + (\phi^{\eta}(\xi,\eta,z))^2 \leq 
            4 \phi^z(\xi,\eta,z)$
	for every  $(\xi,\eta,z) \in V { \times } \erre$;
\item[(c)] $\displaystyle \phi^{\xi} (\xi,\eta, u(\xi,\eta))=2 \diff_{\xi} u (\xi,\eta)$,
	$\displaystyle \phi^{\eta}(\xi,\eta, u(\xi,\eta))=2 \diff_{\eta}
        u (\xi,\eta)$, and
	$\displaystyle \phi^z(\xi,\eta, u(\xi,\eta))= (\diff_{\xi} u (\xi,\eta))^2
	+ (\diff_{\eta}
        u (\xi,\eta))^2$ 
	for every $(\xi,\eta)\in V$;
\item[(d)] $\displaystyle \left( \int_{s}^{t} \phi^{\xi}(\xi,\eta,z)\,
	dz \right)^2 + \left( \int_{s}^{t}
	\phi^{\eta}(\xi,\eta,z)\, dz \right)^2 \leq
	\gamma^2(\xi,\eta)$
	for every $(\xi,\eta) \in V$, $s, t \in \erre$;
\item[(e)] $\displaystyle \int_{u_1}^{u_2} \phi^{\xi}(\xi,0,z)\, dz = 0$ and 
	$\displaystyle \int_{u_1}^{u_2} \phi^{\eta} (\xi,0,z)\, dz = \gamma(\xi,0)=1$
	for every $(\xi,0) \in V$.
\end{description}

Given suitable parameters $\varepsilon >0$ and $\lambda >0$, that
will be chosen later, 
%and assuming
%\begin{equation}
%V = \{ (\xi,\eta): |\eta |< \delta \},
%\end{equation}
we consider the following subsets of $V { \times } \erre$
\begin{eqnarray*}
A_1 & := & \{(\xi,\eta,z)\in V{ \times } \erre : z < u_1(\xi,\eta) - \varepsilon \}, \\
A_2 & := & \{(\xi,\eta ,z)\in V{ \times } \erre : u_1(\xi,\eta) - \varepsilon < z <
	u_1(\xi,\eta) + \varepsilon\}, \\ 
A_3 & := & \{(\xi,\eta,z)\in V{ \times } \erre : u_1(\xi,\eta) + \varepsilon < z <
	\beta_1(\xi,\eta)  \}, \\
A_4 & := & \{(\xi,\eta,z)\in V{ \times } \erre : \beta_1(\xi,\eta)  < z < 
	\beta_2(\xi,\eta)  + 1/\lambda \}, \\
A_5 & := & \{(\xi,\eta,z)\in V{ \times } \erre : \beta_2(\xi,\eta) + 1/\lambda
	< z < u_2(\xi,\eta) - \varepsilon \}, \\
A_6 & := & \{(\xi,\eta ,z)\in V{ \times } \erre : u_2(\xi,\eta) - \varepsilon < z <
	u_2(\xi,\eta) + \varepsilon\}, \\
A_7 & := & \{(\xi,\eta,z)\in V{ \times } \erre : z > u_2(\xi,\eta) + \varepsilon \}, 
\end{eqnarray*}
where $\beta_1$ and $\beta_2$ are suitable smooth function 
such that $u_1(\xi,0)<\beta_1 (\xi,0)= \beta_2(\xi,0)<u_2(\xi,0)$,
which will be defined later. 
Since we suppose $u_2>0$ on $V$, if
$\varepsilon$ is small enough,
while $\lambda$ is sufficiently large, then the sets $A_1, \ldots, A_7$
are nonempty and disjoint, provided $V$ is sufficiently small.
 
The vector 
$\phi(\xi,\eta,z)$
introduced in (\ref{calibra2}) will be written as
$$\phi(\xi,\eta,z)= (\phi^{\xi\eta}(\xi,\eta,z), \phi^z(\xi,\eta,z)),$$
where $\phi^{\xi\eta}$ is the two-dimensional vector
given by the pair $(\phi^{\xi}, \phi^{\eta})$.
For $(\xi,\eta) \in V$ and $z\in\erre$ we define $\phi(\xi,\eta,z)$
as follows:
%
%DEFINIZIONE DEL CAMPO
%
$$
\begin{cases}
\displaystyle (0,\omega_1(\xi,\eta))
& \text{in $A_1 \cup A_3$}, \\
\\
\displaystyle \left( 2 \nablav u_1 
- 2 \frac{u_1 - 
z}{v_1}\nablav v_1,
\left| \nablav u_1 
-  \frac{u_1 - 
z}{v_1}\nablav v_1
\right|^2 \right)
& \text{in $A_2$}, \\
\\
\displaystyle  \left( \lambda \sigma(\xi,\eta) \nablav w , \mu
\right) & \text{in 
$A_4$}, \\
\\
\displaystyle (0, \omega_2(\xi,\eta))
& \text{in $A_5\cup A_7$}, \\
\\
\displaystyle \left( 2 \nablav u_2 
- 2 \frac{u_2 - 
z}{v_2}\nablav v_2,
\left| \nablav u_2 
-  \frac{u_2 - 
z}{v_2}\nablav v_2
\right|^2 \right)
& \text{in $A_6$}, 
\end{cases}$$
where $\nabla$ denotes the gradient with respect to the variables
$(\xi,\eta)$, the functions $v_i$ are defined by 
$$v_1(\xi,\eta):=
\varepsilon +M\eta,\;
v_2(\xi,\eta):=
\varepsilon -M\eta,$$ 
and $M$ and $\mu$ are positive parameters
which will be fixed later, while
\begin{equation}\label{omegai}
\omega_i(\xi,\eta) := 
\frac{\varepsilon^2M^2}{v_i^2(\xi,\eta)} 
- (\partial_{\xi}u_i(\xi,\eta))^2
- (\partial_{\eta}u_i(\xi,\eta))^2
\end{equation}
for $i=1,2$, and for every $(\xi,\eta)\in V$. 
We choose $w$ as the solution of the Cauchy problem
\begin{equation}\label{vudoppio}
\begin{cases}
\triangle w = 0, \\
\displaystyle w(\xi,0) = -\frac{2\varepsilon}{1-2\varepsilon M}\int_0^{\xi}
n(s) (\diff_{\xi} u_1(s,0) + \diff_{\xi} u_2(s,0)) \,ds, \\
\diff_{\eta} w (\xi,0) = n(\xi),
\end{cases}
\end{equation}
where $n$ is a positive analytic function that
will be chosen later in a suitable way (if $V$ is sufficiently small, $w$ is
defined in $V$).
To define $\sigma$, we need some further explanations: we call
$p(\xi,\eta)$ the solution of the problem
\begin{equation}\label{pdef}
\begin{cases}
\displaystyle \diff_{\eta}p(\xi,\eta)=\frac{\diff_{\xi}w}{\diff_{\eta}
w}(p(\xi,\eta),\eta), \\
p(\xi,0)=\xi,
\end{cases}
\end{equation}
which is defined in $V$, provided $V$ is small enough.
By applying the Implicit Function Theorem, it is easy to see
that there exists a function $q$ defined in $V$
(take $V$ smaller, if needed) such that 
\begin{equation}\label{qdef}
p(q(\xi,\eta), \eta)=\xi.
\end{equation}
At last, we define
$$\sigma(\xi,\eta):= \frac{1}{n(q(\xi,\eta))} (1-2\varepsilon M).$$
We choose $\beta_i$, for $i=1,2$, as the solution of the Cauchy problem
\begin{equation}\label{betai}
\begin{cases}
\lambda \sigma(\xi,\eta) \diff_{\xi}w(\xi,\eta) \diff_{\xi}\beta_i(\xi,\eta) +
\lambda \sigma(\xi,\eta) \diff_{\eta}w(\xi,\eta) \diff_{\eta}\beta_i(\xi,\eta)- \mu
= -\omega_i(\xi,\eta), \\
\beta_i(\xi,0)=\displaystyle \frac{1}{2}(u_1(\xi,0)+u_2(\xi,0)).
\end{cases}
\end{equation}
Since the line $\eta=0$ is not characteristic, there exists a unique solution
$\beta_i\in C^{\infty}(V)$, provided $V$ is small enough.

The purpose of the definition of $\phi$ in $A_2$ and $A_6$ is
to provide a divergence free vector field satisfying condition (c)
and such that
\begin{eqnarray*}
& \phi^{\eta}(\xi, 0, z) \geq 0 & \hbox{for}\ u_1<z<u_2, \\
& \phi^{\eta}(\xi, 0, z) \leq 0 & \hbox{for}\ z<u_1 \ \hbox{and}\ z>u_2.
\end{eqnarray*}
These properties are crucial in order to obtain (d) and (e) 
simultaneously.

The role of $A_4$ is to give the main contribution 
to the integral in (e).
The idea of the construction is to start
from the gradient field of a harmonic function $w$ whose
normal derivative is
positive on the line $\eta=0$, while the tangential derivative 
is chosen in order to annihilate
the $\xi$-component of $\phi$, as required in
(e).
Then, we multiply the field by a function $\sigma$ which 
is defined first on $\eta =0$ 
in order to make (e) true, and then
in a neighbourhood of $\eta =0$ by assuming $\sigma$ constant
along the integral curves of the gradient field, so that
$\sigma \nabla w$ remains divergence free.

The other sets $A_i$ are simply regions of transition, where the field is 
taken purely vertical.

Let us prove condition (a).
By Lemma \ref{alberti}
it follows that $\phi$ is divergence free in 
$A_2 \cup A_6$, noting that it is constructed starting from
the family of harmonic functions $u_i(\xi,\eta) - t 
v_i(\xi, \eta )$.

In $A_4$ condition (a) is true since, as remarked above, $\phi$ is the product
of $\nabla w$ with the function $\sigma$ which, by construction, is constant along
the integral curves of $\nabla w$.

In the other sets, condition (a) is trivially satisfied.

Note that the normal component of $\phi$ is continuous across each $\partial A_i$:
for the regions $A_2$, $A_6$,  and 
for $A_4$, this continuity is guaranteed by our choice of $\omega_i$ and $\beta_i$
respectively. This implies
that (a) is satisfied in the sense of distributions on $V{ \times}\erre$.

Since 
$$\omega_i(\xi,0) = M^2- (\diff_{\xi}u_i(\xi,0))^2,$$
condition (b) is satisfied in $A_1\cup A_3$ and in $ A_5\cup A_7$ if we require
that
$$M>\sup \{ |\diff_{\xi}u_i(\xi,0)| :\ (\xi,0)\in V,\ i=1,2 \},$$
provided $V$ is small enough.

Arguing in a similar way, if we impose that
$$\mu  >  \sup \left\{ \frac{\lambda^2}{4}(1-2 \varepsilon M)^2 \left(1 +
\frac{4\varepsilon^2}{(1-2\varepsilon M)^2}(\diff_{\xi}u_1(\xi,0)+
\diff_{\xi}u_2(\xi,0))^2 \right) : \ (\xi,0)\in V \right\},$$
condition (b) holds in $A_4$, provided $V$ is sufficiently small.

In the other cases, (b) is trivial. 

Looking at the definition of $\phi$ on $A_2$ and $A_6$, one can check
that condition (c) is satisfied.

By direct computations we find that 
\begin{eqnarray}\label{fixi}
\int_{u_1}^{u_2}\phi^{\xi}\ dz & = &  
2\varepsilon \diff_{\xi} u_1 + 
2\varepsilon \diff_{\xi} u_2 + 
\lambda
\left(\beta_2-\beta_1+\frac{1}{\lambda}\right)\sigma
\diff_{\xi}w,\\
\int_{u_1}^{u_2}\phi^{\eta} \label{fieta}
\ dz & = &
2 \varepsilon \diff_{\eta}u_1 + 
2 \varepsilon \diff_{\eta} u_2 + 
M \frac{\varepsilon^2}{\varepsilon + M\eta}
+ M \frac{\varepsilon^2}{\varepsilon - M\eta}
+\lambda
\left(\beta_2-\beta_1+\frac{1}{\lambda}\right)\sigma
\diff_{\eta}w,
\end{eqnarray}
for every $(\xi,\eta)\in V$.

By using (\ref{vudoppio})
and the definition of $\sigma$, we obtain
\begin{equation}\label{e/2}
\int_{u_1(\xi,0)}^{u_2(\xi,0)}\phi^{\xi}
(\xi,0,z)\ dz=0
\end{equation}
and
\begin{equation}\label{conde}
\int_{u_1(\xi,0)}^{u_2(\xi,0)}\phi^{\eta}
(\xi,0,z)\ dz=1,
\end{equation} 
so condition (e) is satisfied.

The proof of condition (d) will be split in two steps:
we first prove
that condition (d) holds if $s$ and $t$ respectively belong
to a suitable neighbourhood of $u_1(\xi,\eta)$ and $u_2(\xi,\eta)$, 
whose width is uniform
with respect to $(\xi,\eta)$ in $V$; then,
by a quite simple continuity argument we show 
that condition (d) is true 
if $s$ or $t$ is not too close 
to $u_1(\xi, \eta)$ or $u_2(\xi,\eta)$ respectively.

For $(\xi, \eta)\in V$ and $s, t\in \erre$,
we set
$$\I{\eta}{s}{t} :=
\int_{s}^{t} \phi^{\xi\eta}(\xi, \eta, z)\,dz$$
and we denote by $I^{\xi}$ and $I^{\eta}$ its components.

\bstep{1}{For a suitable choice of $\varepsilon$ and of the function $n$
(see (\ref{vudoppio}))
there exists $\delta >0$ such that condition (d) holds
for $|s-u_1(\xi,\eta)|<\delta$, $|t-u_2(\xi,\eta)|<\delta$, and
$(\xi,\eta)\in V$, provided $V$ is small enough.}

\noindent To estimate the vector whose components are given by (\ref{fixi}) and (\ref{fieta}),
we use suitable polar coordinates. If $V$ is small enough, for every $(\xi,\eta)\in V$
there exist $\rho\en(\xi,\eta)>0$ and $-\pi/2<\theta\en(\xi,\eta)<\pi/2$ such that
\begin{eqnarray}
I^{\xi}(\xi,\eta,u_1(\xi,\eta),u_2(\xi,\eta)) & = & 
\emme(\xi,\eta)\sin \theta\en(\xi,\eta), \label{pol1} \\
I^{\eta}(\xi,\eta,u_1(\xi,\eta),u_2(\xi,\eta)) & = &
\emme(\xi,\eta)\cos \theta\en(\xi,\eta).\label{pol2}
\end{eqnarray}
In the notation above we have made explicit the dependence on the parameter $\varepsilon$
and on the function $n$ which appears in the definition of $w$ (see (\ref{vudoppio})).

In order to prove condition (d), we want to compare the behaviour of the functions $\emme$ 
and $\gamma$ for $|\eta|$ small.
We have already proved that $\emme(\xi,0)=\gamma(\xi,0)=1$; 
we start computing the first
derivative of $\gamma$ and of $\emme$ with respect to the variable $\eta$.

\bclaim{1}{$\diff_{\eta}(|\nabla_{\!xy}\xi(\Psi)|^2)(\xi,0)=-2\,
\curv\Gamma(\xi)$.}

\eclaim{ By (\ref{gi}) we obtain
$$|\nabla_{\!xy}\xi(\Psi)|^2=\frac{1}{(\diff_{\xi}\varx)^2+(\diff_{\xi}\vary)^2},$$
hence
\begin{equation}\label{002}
\diff_{\eta}(|\nabla_{\!xy}\xi(\Psi)|)^2=-[(\diff_{\xi}\varx)^2+(\diff_{\xi}\vary)^2]
^{-2}(2\diff_{\xi}\varx\,\diff^2_{\xi\eta}\varx+2\diff_{\xi}\vary\,
\diff^2_{\xi\eta}\vary).
\end{equation}
Using the fact that $(\diff_{\xi}\varx)^2+(\diff_{\xi}\vary)^2$ is equal to $1$ 
at $(\xi,0)$, and the equalities in (\ref{001}), we finally get
$$\diff_{\eta}(|\nabla_{\!xy}\xi(\Psi)|^2)(\xi,0)=-2(-\diff_{\xi}\varx\,
\diff^2_{\xi\xi}\vary
+\diff_{\xi}\vary\,\diff^2_{\xi\xi}\varx)=-2\,\curv\Gamma(\xi),$$
where the last equality follows from (\ref{curvatura1}): therefore the claim
is proved.}
 
Since $\gamma =(|\nabla_{\!xy}\xi(\Psi)|^2)^{-\frac{1}{2}}$, one
has that $\diff_{\eta}\gamma =
-\frac{1}{2}
(|\nabla_{\!xy}\xi(\Psi)|^2)^{-\frac{3}{2}}\diff_{\eta}(|\nabla_{\!xy}\xi(\Psi)|^2)$; 
using the previous claim we can conclude that  $$\diff_{\eta}(\gamma)
(\xi,0)=-\frac{1}{2}\diff_{\eta}(|\nabla_{\!xy}\xi(\Psi)|^2)(\xi,0)=
\curv\Gamma(\xi) .$$
Using the equality
$$\emme^2(\xi, \eta) =\left(  
I^{\xi}(\xi,\eta,u_1(\xi,\eta),u_2(\xi,\eta))
\right)^2 + \left( 
I^{\eta}(\xi,\eta,u_1(\xi,\eta),u_2(\xi,\eta))
\right)^2,$$
we obtain
\begin{equation}\nonumber
\diff_{\eta}(\emme)  =
\frac{1}{\emme} \diff_{\eta}\left(
I^{\xi}(\xi,\eta,u_1,u_2)
\right) I^{\xi}(\xi,\eta,u_1,u_2) 
+ \frac{1}{\emme} \diff_{\eta} \left( 
I^{\eta}(\xi,\eta,u_1,u_2)
\right)I^{\eta}(\xi,\eta,u_1,u_2).
\end{equation}
By (\ref{e/2}) it follows that the first addend in the  
expression above is equal to zero at $(\xi,0)$, while by (\ref{conde})
it turns out that $I^{\eta}(\xi,0,u_1,u_2)= \emme(\xi,0) =1$;
therefore,
\begin{equation}\label{2000}
\diff_{\eta}(\emme)(\xi,0)  =
\diff_{\eta} \left( 
I^{\eta}(\xi,0,u_1,u_2)
\right).
\end{equation}
By (\ref{fieta}) it follows that
\begin{eqnarray}
\diff_{\eta} \left( 
I^{\eta}(\xi,\eta,u_1,u_2)
\right) \hspace{-.3cm} & =& \hspace{-.3cm} 
2\varepsilon \diff^2_{\eta} u_1 + 2\varepsilon \diff^2_{\eta} u_2
-\frac{\varepsilon^2}{(\varepsilon +M\eta)^2}M^2
+\frac{\varepsilon^2}{(\varepsilon - M\eta)^2}M^2 +
\lambda(\diff_{\eta}\beta_2- \diff_{\eta}\beta_1)\sigma \diff_{\eta}w +
\nonumber \\ 
& & \hspace{-.3cm} +\lambda(\beta_2-\beta_1+1/\lambda)\diff_{\eta}
(\sigma\diff_{\eta}w). \label{400}
\end{eqnarray}
From (\ref{betai}) and the
Euler condition iii),
we have that
\begin{eqnarray}\label{dbetai}
\lambda(\diff_{\eta}\beta_2(\xi,0)-\diff_{\eta}\beta_1(\xi,0))\sigma(\xi,0)
\diff_{\eta}w(\xi,0) & = &  -\omega_2(\xi,0)
+\omega_1(\xi,0) \nonumber \\
& = &
(\diff_{\xi}u_2(\xi,0))^2-(\diff_{\xi}u_1(\xi,0))^2 \nonumber \\ 
& = & \curv\Gamma(\xi),
\end{eqnarray}
while 
$$\diff_{\eta}(\sigma \diff_{\eta}w)(\xi,0) =
- \diff_{\xi}(\sigma \diff_{\xi}w) (\xi,0) =
 \diff_{\xi} (2 \varepsilon \diff_{\xi} u_1 (\xi,0) + 
2 \varepsilon \diff_{\xi} u_2)(\xi,0),$$
where we have used the fact that $\sigma\nabla w$ is divergence free and
the definition of $\sigma$ and $w$.
Putting this last fact together with (\ref{400}), (\ref{dbetai}), and the
harmonicity of $u_i$, we finally get
\begin{equation}\label{uguali}
\diff_{\eta}(\emme)(\xi,0)= 
\curv\Gamma(\xi) =\diff_{\eta}(\gamma)(\xi,0). 
\end{equation}

\bclaim{2}{ $\diff^2_{\eta\eta}(|\nabla_{\!xy}\xi(\Psi)|^2)(\xi,0)=
4\,[\curv\Gamma(\xi)]^2.$}

\eclaim{By differentiating with respect to $\eta$ the expression in
(\ref{002}) and by (\ref{001}), we obtain 
\begin{eqnarray*}
\diff^2_{\eta\eta}(|\nabla_{\!xy}\xi(\Psi)|^2) 
& = &  - 2 [(\diff_{\xi}\varx)^2
+(\diff_{\xi}\vary)^2]^{-2}[(\diff_{\xi\eta}\varx)^2+\diff_{\xi}\varx\,
\diff^3_{\xi\eta\eta}\varx+(\diff_{\xi\eta}\vary)^2+\diff_{\xi}\vary\,
\diff^3_{\xi\eta\eta}\vary] + \\
& & + 8[(\diff_{\xi}\varx)^2
+(\diff_{\xi}\vary)^2]^{-3}
(\diff_{\xi}\varx\,\diff^2_{\xi\eta}\varx +
\diff_{\xi}\vary\,\diff^2_{\xi\eta}\vary )^2 \\
& = &  -2 [(\diff_{\xi}\varx)^2
+(\diff_{\xi}\vary)^2]^{-2}[(\diff_{\xi\xi}^2\vary)^2+(\diff_{\xi\xi}^2\varx)^2 
-\diff_{\xi}\varx\,\diff^3_{\xi\xi\xi}\varx-
\diff_{\xi}\vary\,\diff^3_{\xi\xi\xi}\vary] + \\
& & + 8[(\diff_{\xi}\varx)^2
+(\diff_{\xi}\vary)^2]^{-3}
(- \diff_{\xi}\varx\,\diff^2_{\xi\xi}\vary +
\diff_{\xi}\vary\,\diff^2_{\xi\xi}\varx )^2.
\end{eqnarray*}
Note that
$$-\diff_{\xi}\varx\,\diff^3_{\xi\xi\xi}\varx-
\diff_{\xi}\vary\,\diff^3_{\xi\xi\xi}\vary =
(\diff_{\xi\xi}^2\vary)^2 + (\diff_{\xi\xi}^2\varx)^2 
-\frac{1}{2}
\diff^2_{\xi\xi}((\diff_{\xi}\varx)^2+(\diff_{\xi}\vary)^2).$$
Using (\ref{curvatura1}), (\ref{curvq}), and the fact that
$(\diff_{\xi}\varx)^2+(\diff_{\xi}\vary)^2$ is equal to $1$ at $(\xi,0)$,
we obtain the claim.}

By using Claims~1 and 2, we can conclude that
\begin{eqnarray}
\diff^2_{\eta\eta}(\gamma)(\xi,0) & = & \left.
\left[\frac{3}{4}(|\nabla_{\!xy}\xi(\Psi)|^2)^{-\frac{5}{2}}
[\diff_{\eta}(|\nabla_{\!xy}\xi(\Psi)|^2)]^2
- \frac{1}{2}(|\nabla_{\!xy}\xi(\Psi)|^2)^{-\frac{3}{2}}
\diff^2_{\eta\eta}(|\nabla_{\!xy}\xi(\Psi)|^2)
\right]\right|_{(\xi,0)} \nonumber\\
&=& [\curv\Gamma(\xi)]^2.\label{101}
\end{eqnarray}
The second derivative of $\emme$ with respect to $\eta$ is given by
\begin{eqnarray*}
\diff^2_{\eta\eta}(\emme) & = &
\frac{1}{\emme} \left\{ \left[ \diff_{\eta} \left(
I^{\xi}(\xi,\eta,u_1,u_2)
\right) \right]^2 +\diff^2_{\eta\eta} \left(
I^{\xi}(\xi,\eta,u_1,u_2)
\right) I^{\xi}(\xi,\eta,u_1,u_2) 
 + \right. \\
& & \left.  + \left[ \diff_{\eta} \left(
I^{\eta}(\xi,\eta,u_1,u_2)
\right) \right]^2 
+ \diff^2_{\eta\eta} \left(
I^{\eta}(\xi,\eta,u_1,u_2)
\right) I^{\eta}(\xi,\eta,u_1,u_2) \right\} - \\
& &  - \frac{1}{\emme} [\diff_{\eta}(\emme)]^2. 
\end{eqnarray*}
By the equalities (\ref{e/2}), (\ref{conde}), and (\ref{2000}),
the expression above computed at $(\xi,0)$ reduces to 
\begin{equation}\label{404}
\diff^2_{\eta\eta}(\emme)(\xi,0) =
\left[ \diff_{\eta} \left( \left.
I^{\xi}(\xi,\eta,u_1,u_2)
\right)\right|_{(\xi,0)} \right]^2
+ \diff^2_{\eta\eta} \left( \left.
I^{\eta}(\xi,\eta,u_1,u_2)
\right)\right|_{(\xi,0)}. 
\end{equation}
By differentiating (\ref{fixi}) and (\ref{400})
with respect to $\eta$, 
we obtain that 
\begin{equation}\label{402}
\diff_{\eta} \left(
I^{\xi}(\xi,\eta,u_1,u_2)
\right)(\xi,0) = [ 
\lambda(\diff_{\eta} \beta_2 -\diff_{\eta} \beta_1)\sigma \diff_{\xi} w
+ \diff_{\eta}\sigma\,
\diff_{\xi}w + \sigma \diff^2_{\xi\eta} w
] |_{(\xi,0)},
\end{equation}
and
\begin{eqnarray}
\diff^2_{\eta\eta} \left(
I^{\eta}(\xi,\eta,u_1,u_2)
\right)(\xi,0) & = & \frac{4}{\varepsilon} M^3
+ \lambda
[\diff^2_{\eta\eta}\beta_2(\xi,0)-
\diff^2_{\eta\eta}\beta_1(\xi,0)]\sigma(\xi,0)\diff_{\eta}w(\xi,0) + \nonumber \\
& &
+ 2\lambda [\diff_{\eta}\beta_2(\xi,0)-
\diff_{\eta}\beta_1(\xi,0)]\diff_{\eta}( \sigma\diff_{\eta}w)(\xi,0)
+  \diff^2_{\eta\eta}\sigma(\xi,0) \diff_{\eta} w(\xi,0)+ \nonumber \\
& &+ 2 \diff_{\eta}\sigma(\xi,0)\diff^2_{\eta\eta} w(\xi,0)
+ \sigma(\xi,0)\diff^3_{\eta\eta\eta}w(\xi,0),\label{403}  
\end{eqnarray}
while, by using the equation (\ref{betai}),
\begin{eqnarray*}
[\lambda (\diff^2_{\eta\eta}\beta_2-\diff^2_{\eta\eta}\beta_1)
\sigma\diff_{\eta}w]|_{(\xi,0)} & =\hspace{-.3cm}  &
[\diff_{\eta}\omega_1 - \diff_{\eta} \omega_2
- \lambda \diff_{\eta}(\diff_{\xi} \beta_2 - \diff_{\xi} \beta_1)\sigma
\diff_{\xi}w
- \lambda \diff_{\eta}( \sigma \diff_{\eta}w)(\diff_{\eta} \beta_2 - \diff_{\eta} \beta_1)]
|_{(\xi,0)} \\
& =\hspace{-.3cm}  & [-\frac{4}{\varepsilon}M^3 -
\lambda \diff_{\xi}(\diff_{\eta} \beta_2 - \diff_{\eta} \beta_1)\sigma
\diff_{\xi}w
+ \lambda \diff_{\xi}( \sigma \diff_{\xi}w)(\diff_{\eta} \beta_2 - \diff_{\eta} \beta_1)]
|_{(\xi,0)}.
\end{eqnarray*}
Since by (\ref{dbetai}) and by the definition of $\sigma$ we have that
$$\lambda [\diff_{\eta} \beta_2(\xi,0) - \diff_{\eta} \beta_1(\xi,0)]
= \frac{\curv \Gamma(\xi)}{1-2\varepsilon M},$$
and moreover,
$$\sigma(\xi,0) \diff_{\xi}w(\xi,0) = -2\varepsilon
(\diff_{\xi}u_1(\xi,0) + \diff_{\xi}u_2(\xi,0)),$$
we obtain that
\begin{multline}\nonumber
[\lambda (\diff^2_{\eta\eta}\beta_2-\diff^2_{\eta\eta}\beta_1)
\sigma\diff_{\eta}w
+ 2\lambda (\diff_{\eta}\beta_2 -
\diff_{\eta}\beta_1) \diff_{\eta}( \sigma\diff_{\eta}w)]|_{(\xi,0)}
= \\
=-\frac{4}{\varepsilon}M^3
+ \frac{2\varepsilon}{1-2\varepsilon M} 
\diff_{\xi}((\diff_{\xi}u_1-\diff_{\xi}u_2)\,\curv \Gamma )(\xi,0).
\end{multline}
By using the definition of $\sigma$,
we can write
\begin{eqnarray*}
\diff_{\eta} \sigma & = &
- (1 - 2\varepsilon M) \frac{n'(\xi)}{n^2(\xi)} \diff_{\eta}q, \\
\diff^2_{\eta\eta} \sigma & = &
- (1 - 2\varepsilon M) \left[ - 2 \frac{(n'(\xi))^2}{n^3(\xi)} (\diff_{\eta}q)^2
+ \frac{n''(\xi)}{n^2(\xi)} (\diff_{\eta}q)^2
+ \frac{n'(\xi)}{n^2(\xi)} \diff^2_{\eta}q \right].
\end{eqnarray*}
In order to compute the derivatives of $q$, we differentiate the
equality (\ref{qdef}) with respect to $\eta$:
\begin{eqnarray*}
\diff_{\eta} q(\xi,0) & = & - \diff_{\eta} p(\xi,0)
= \frac{2\varepsilon}{1-2\varepsilon M} (\diff_{\xi} u_1(\xi,0)
+ \diff_{\xi} u_2(\xi,0)), \\
\diff^2_{\eta\eta} q(\xi,0) & = & - 2 \diff^2_{\xi\eta} p(\xi,0)
\diff_{\eta} q(\xi,0) - \diff^2_{\eta\eta} p(\xi,0)  \\
& = & \left[-\frac{(\diff_{\xi}w)^2}{(\diff_{\eta}w)^3}\diff^2_{\xi\eta}w
-\frac{1}{\diff_{\eta}w}\diff^2_{\xi\eta}w \right](\xi,0).
\end{eqnarray*}
By the definition of $w$, we obtain
$$\diff^2_{\eta} q(\xi,0) = -\frac{n'(\xi)}{n(\xi)} - \frac{n'(\xi)}{n(\xi)}
\frac{4\varepsilon^2}{(1-2\varepsilon M)^2}
(\diff_{\xi}u_1(\xi,0)+\diff_{\xi}u_2(\xi,0))^2.$$
Finally, we have
\begin{eqnarray*}
\diff^2_{\eta\eta}w (\xi,0) & = &
- \diff^2_{\xi\xi}w(\xi,0) =
\frac{2\varepsilon}{1-2\varepsilon M} [n'(\diff_{\xi}u_1 + \diff_{\xi}u_2)
+ n (\diff^2_{\xi\xi}u_1 + \diff^2_{\xi\xi}u_2)]|_{(\xi,0)}, \\
\diff^3_{\eta\eta\eta}w(\xi,0) & = &
- \diff^2_{\xi\xi}\diff_{\eta}w(\xi,0)
= - n''(\xi).
\end{eqnarray*}
By substituting all information above in (\ref{402}) and in (\ref{403}),
and by using (\ref{404}), we finally obtain that
\begin{eqnarray}
\diff^2_{\eta\eta}(\emme)(\xi,0)
& = & - a_{\varepsilon}(\xi) \frac{n''(\xi)}{n(\xi)}
	+ h_{\varepsilon}\left(\xi, \frac{n'(\xi)}{n(\xi)} \right)\nonumber \\
& = & - a_{\varepsilon}(\xi) \left( \frac{n'(\xi)}{n(\xi)} \right)'
 + h_{\varepsilon}\left(\xi, \frac{n'(\xi)}{n(\xi)} \right) - a_{\varepsilon}(\xi)
 \left(\frac{n'(\xi)}{n(\xi)}\right)^2, \label{3000}
\end{eqnarray}
where 
\begin{equation}
\begin{array}{cl}
 a_{\varepsilon} (\xi) \to 1  & \hbox{uniformly in $[0,l(\Gamma)]$,}
\vspace{.1cm}  \\
  h_{\varepsilon} (\xi, \tau) \to 2\tau^2 & \hbox{uniformly on 
the compact sets of $[0, l(\Gamma)]{ \times } \erre$,} \label{ahconv}
\end{array}
\end{equation}
as $\varepsilon \to 0$.

\bclaim{3}{There exists $\overline{\varepsilon}>0$
such that for every $\varepsilon\in (0,\overline{\varepsilon})$,
we can find
an analytic function $n : [0, l(\Gamma)] \to (0,+\infty)$ satisfying
\begin{equation}\label{diverse}
\diff^2_{\eta\eta}
(\emme -\gamma)(\xi,0)=-\frac{\pi^2}{16\, l^2(\Gamma)}
\qquad \hbox{and} \qquad
\left|\frac{n'(\xi)}{n(\xi)}\right| \leq N
\quad \forall \xi\in [0, l(\Gamma)],
\end{equation}
where $N:=1+\max \left\{ \displaystyle \frac{\pi}{4\, l(\Gamma)} ,
k(\Gamma ) \right\}$ and $k(\Gamma)=\| \curv\Gamma \|_{\infty}$.}

\eclaim{Set $\tau:=n'/n$; in order to prove the claim, by (\ref{3000}) 
and (\ref{101}) we
study the Cauchy problem
\begin{equation}\label{eqvera}
\begin{cases}
-a_{\varepsilon}(\xi)\tau'+
h_{\varepsilon}(\xi,\tau) -\tau^2
-[\curv\Gamma(\xi)]^2=\displaystyle -\frac{\pi^2}{16\, l^2(\Gamma)}, \\
\tau(0)=0,
\end{cases}
\end{equation}
and we investigate for which values of $\varepsilon$
it admits a solution defined in the whole interval $[0, l(\Gamma)]$, with
$L^{\infty}$-norm less than $N$. 
As $\varepsilon\to 0$, by (\ref{ahconv})
we obtain the limit problem
\begin{equation}\label{eqlim}
\begin{cases}
-\tau'+
\tau^2 -(\curv\Gamma)^2 =\displaystyle -\frac{\pi^2}{16\, l^2(\Gamma)} ,\\
\tau(0)=0.
\end{cases}
\end{equation}
By comparing with the solutions $\tau_1$ and $\tau_2$ of the Cauchy problems
\begin{equation}\label{auxil}
\begin{cases}
-\tau'_1+ \tau_1^2 = \displaystyle -\frac{\pi^2}{16\, l^2(\Gamma)}, \\
\tau_1(0)=0,
\end{cases} 
\qquad  
\begin{cases}
-\tau'_2+
\tau_2^2 - k^2(\Gamma) = 
\displaystyle -\frac{\pi^2}{16\, l^2(\Gamma)}, \\
\tau_2(0)=0,
\end{cases} 
\end{equation}
one easily sees that the
solution of (\ref{eqlim}) is defined in $[0, l(\Gamma)]$, 
with $L^{\infty}$-norm less than the maximum between $\| \tau_1 \|_{\infty}$
and $\| \tau_2 \|_{\infty}$, which is, by explicit computation, 
less than $\max \{\pi/(4 l(\Gamma)) ,
 k(\Gamma) \}$.
By the theorem of continuous dependence on the coefficients (see
\cite{Har}), we can find $\overline{\varepsilon}$
such that, for every $\varepsilon\in(0,\overline{\varepsilon})$, the
solution of (\ref{eqvera}) is defined in
$[0, l(\Gamma)]$ with $L^{\infty}$-norm less than $N$.}

For every $\varepsilon\in (0, \overline{\varepsilon})$,
we set
\begin{equation}\label{neps}
n_{\varepsilon}(\xi):= {\rm e}^{\int_0^{\xi} \tau_{\varepsilon}(s)\,ds},
\end{equation}
where $\tau_{\varepsilon}$ is the solution of (\ref{eqvera}).

From now on we will simply write $\emmee$ and $\thetae$ instead of
$\rho_{\varepsilon, n_{\varepsilon}}$
and $\thetak$.

We now want to estimate the angle $\thetae(\xi,\eta)$ by a quantity which is
independent of $\varepsilon$. Since by (\ref{fixi}) and (\ref{fieta})
$$\tan\thetae =\frac{2\varepsilon \diff_{\xi} u_1 + 
2\varepsilon \diff_{\xi} u_2 + 
\lambda
\left(\beta_2-\beta_1+\frac{1}{\lambda}\right)\sigma
\diff_{\xi}w}
{2 \varepsilon \diff_{\eta}u_1 + 
2 \varepsilon \diff_{\eta} u_2 + 
M \varepsilon^2 (\varepsilon + M\eta)^{-1}
+ M \varepsilon^2 (\varepsilon - M\eta)^{-1}
+\lambda
\left(\beta_2-\beta_1+\frac{2}{\lambda}\right)\sigma
\diff_{\eta}w},$$
we have
$$\diff_{\eta}\thetae(\xi,0)=
-\frac{2\varepsilon}{1-2\varepsilon M}(\diff_{\xi}
u_1+\diff_{\xi}u_2)
\left( \curv \Gamma - 
2\varepsilon (\diff_{\xi}
u_1 +\diff_{\xi}u_2 )
\frac{n'_{\varepsilon}(\xi)}{n_{\varepsilon}(\xi)} \right)
+ (1-2\varepsilon M)\frac{n'_{\varepsilon}(\xi)}{n_{\varepsilon}(\xi)},$$
and so, by Claim~3, if $\varepsilon$ is sufficiently small,
\begin{equation}\label{deranglim}
|\diff_{\eta}\thetae(\xi,0)| < N \qquad \forall\xi\in [0, l(\Gamma)].
\end{equation}
Let $\thetat(\eta)$ be an arbitrary continuous function with 
\begin{equation}\label{thtilde}
\thetat(0)=0 \qquad \hbox{and}  \qquad \thetat'(0) = N;
\end{equation}
by (\ref{deranglim}), it follows that
\begin{equation}\label{angolo}
|\thetae(\xi,\eta)|<\thetat(\eta)\,\hbox{sign}\,\eta
\end{equation}
for every $(\xi,\eta)\in V$, provided $V$ is sufficiently small.

Given $h>0$, we consider the vectors
\begin{eqnarray*}
b^h_1(\xi, \eta, s) & := & \left( 0,-2
(s-u_1(\xi,\eta))\diff_{\eta}u_1(\xi,\eta)  -h(s-u_1(\xi,\eta))^2 \right), \\
b^h_2(\xi, \eta, t) & := & \left( 0,2
(t-u_2(\xi,\eta))\diff_{\eta}u_2(\xi,\eta)  -h(t-u_2(\xi,\eta))^2 \right)
\end{eqnarray*}
for $(\xi, \eta)\in V$ and $s, t\in \erre$. We denote by $B(r)$ the open ball centred at $(0,-r)$
with radius $r$.

Let us define 
$\ren(\xi,\eta,s,t)$ as the maximum radius $r$ such that the set
$$(\emmee(\xi,\eta)\sin\thetat(\eta),\emmee
(\xi,\eta)\cos\thetat(\eta))+b_1^h(\xi,\eta,s)+b_2^h(\xi,\eta,t) + B(r)$$
is contained in the ball centred at $(0,0)$ with radius $\gamma(\xi,\eta)$. 

\bclaim{4}{If we define
\begin{equation}\label{defd}
d:= \frac{1}{1+ 16 \, l^2(\Gamma) N^2/\pi^2},
\end{equation}
where $N$ is the constant introduced in the previous claim,
then there exists $h>0$ such
that for every $\varepsilon\in (0,\overline{\varepsilon})$ (see Claim~3), there
exists $\delta\in(0,\varepsilon)$ so that,
if $V$ is small enough,
\begin{equation}\label{infpos}
\inf\left\{2\,\ren(\xi,\eta,s,t) :
(\xi,\eta)\in V,\ |s-u_1(\xi,\eta)|\leq\delta,\,
|t-u_2(\xi,\eta)|\leq\delta\right\} > \frac{d}{2}.
\end{equation}}

\eclaim{Let $\overline{\rho}_{\varepsilon}^h(\xi,\eta,s,t)>0$ and 
$ -\pi/2<\stheta(\xi,\eta,s,t)<\pi/2$ be such that
\begin{multline}\label{defover}
\left(\emmee(\xi,\eta)\sin\thetat(\eta),\emmee
(\xi,\eta)\cos\thetat(\eta)\right)+b_1^h(\xi,\eta,s)+ b_2^h(\xi,\eta,t) = \\
=\left(\overline{\rho}_{\varepsilon}^h(\xi,\eta,s,t)\sin\stheta(\xi,\eta,s,t),
\overline{\rho}_ {\varepsilon}^h(\xi,\eta,s,t)\cos\stheta(\xi,\eta,s,t)\right).
\end{multline}
To prove Claim~4, it is enough to show that, for every
$\varepsilon\in(0,\overline{\varepsilon})$,   
there exists $\delta\in(0,\varepsilon)$ with the property that
\begin{equation}\label{disug}
\left(1-
\frac{d}{2}\cos\stheta(\xi,\eta,s,t)\right)\overline{\rho}_{\varepsilon}^h(\xi,\eta,s,t)<  
\left( 1-\frac{d}{2}\right)\gamma(\xi,\eta)
\end{equation} 
for $|s-u_1(\xi,\eta)|\leq\delta$,
$|t-u_2(\xi,\eta)|\leq\delta$, and $(\xi,\eta)\in V$ with $\eta\neq0$,
provided $V$ is sufficiently small. Indeed, if (\ref{disug}) holds, it follows
in particular that $\overline{\rho}_{\varepsilon}^h(\xi,\eta,s,t)<
\gamma(\xi,\eta)$, and this inequality with some easy geometric computations
implies that 
$$2\,\ren(\xi,\eta,s,t)= \frac{\gamma^2(\xi,\eta)-
(\overline{\rho}_{\varepsilon}^h(\xi,\eta,s,t))^2}{\gamma-
\overline{\rho}_{\varepsilon}^h(\xi,\eta,s,t)\cos\stheta(\xi,\eta,s,t)};$$ at
this point, it is easy to see that, if $V$ is small enough,  
inequality (\ref{disug}) implies that $2\,\ren(\xi,\eta,s,t)>d/2$, that is
Claim~4. So let us prove (\ref{disug}).

We set 
$$f^{d,h}(\xi,\eta,s,t): =
\left(1-\frac{d}{2}\cos\stheta(\xi,\eta,s,t)\right)\overline{\rho}_{\varepsilon}^h(\xi,\eta,s,t)- 
\left( 1-\frac{d}{2}\right)\gamma(\xi,\eta)$$
and we note that $f^{d,h}(\xi,0,u_1(\xi,0),u_2(\xi,0))=0$. 
We will show that 
\begin{enumerate}
\item $\nabla_{\!\eta st}f^{d,h}(\xi,0,u_1(\xi,0),u_2(\xi,0))=0$ if $(\xi,0)\in V$,
\item $\nabla^2_{\!\eta st}f^{d,h}(\xi,0,u_1(\xi,0),u_2(\xi,0))$ 
is negative definite if $(\xi,0)\in V$,
\end{enumerate}
where $\nabla_{\!\eta st}f^{d,h}$ and $\nabla^2_{\!\eta st}f^{d,h}$
denote respectively the gradient and the hessian matrix of $f^{d,h}$
with respect to the variables $(\eta,s,t)$.
Equality $1$ follows by direct computations and by (\ref{uguali}).
Using (\ref{defover}), the equality in (\ref{diverse}), 
and (\ref{thtilde}), we obtain
$$\diff^2_{\eta\eta}f^{d,h}(\xi,0,u_1(\xi,0),u_2(\xi,0))
=-\frac{\pi^2}{16\, l^2(\Gamma)}\left(1-
\frac{d}{2}\right) +
\frac{d}{2}N^2;$$ 
then by the definition of $d$,
\begin{equation}\label{primo}
\diff^2_{\eta\eta}f^{d,h}(\xi,0,u_1(\xi,0),u_2(\xi,0))
= - \frac{\pi^2}{32\, l^2(\Gamma)}
<0.
\end{equation}
Moreover we easily obtain that 
$$\diff^2_{tt}f^{d,h}(\xi,0,u_1(\xi,0),u_2(\xi,0))=
\diff^2_{ss}f^{d,h}(\xi,0,u_1(\xi,0),u_2(\xi,0))=-2h\left( 1- \frac{d}{2}
\right),$$
$$\diff^2_{s\eta}f^{d,h}(\xi,0,u_1(\xi,0),u_2(\xi,0))=
-2\left( 1- \frac{d}{2}
\right)\diff^2_{\eta\eta}u_1(\xi,0),$$ 
$$\diff^2_{t\eta}f^{d,h}(\xi,0,u_1(\xi,0),u_2(\xi,0))=2\left( 1- \frac{d}{2}
\right)
\diff^2_{\eta\eta}u_2(\xi,0),$$
$$\diff^2_{ts}f^{d,h}(\xi,0,u_1(\xi,0),u_2(\xi,0))=0.$$
By the expressions, it follows that
\begin{equation}\nonumber
\det\left( 
\begin{array}{cc}
\diff^2_{\eta\eta}f^{d,h} & \diff^2_{s\eta}f^{d,h} \\
\diff^2_{s\eta}f^{d,h} & \diff^2_{ss}f^{d,h}
\end{array}
\right)(\xi,0,u_1(\xi,0),u_2(\xi,0))= h ( 2- d)
\frac{\pi^2}{32\, l^2(\Gamma)}
-( 2- d )^2[\diff^2_{\eta\eta}u_1(\xi,0)]^2,
\end{equation}
and that the determinant of the hessian matrix of $f^{d,h}$ at
$(\xi,0,u_1(\xi,0),u_2(\xi,0))$ is given by
$$\det \nabla^2_{\!\eta st}f^{d,h}(\xi,0,u_1(\xi,0),u_2(\xi,0))=
-h^2( 2- d)^2
\frac{\pi^2}{32\, l^2(\Gamma)}
+h( 2- d)^3[(\diff^2_{\eta\eta}u_1(\xi,0))^2
+ (\diff^2_{\eta\eta}u_2(\xi,0))^2].$$
By the definition of $d$, if $h$ satisfies
\begin{equation}\label{noname}
h> \frac{32}{\pi^2} (2- d) l^2(\Gamma)
\sum_{i=1}^{2}\| \diff^2_{\eta\eta}u_i \|_{L^{\infty}(\Gamma)}^2,
\end{equation}
then for every $(\xi,0)\in V$ we have
\begin{equation}\label{secondo}
\det\left( 
\begin{array}{cc}
\diff^2_{\eta\eta}f^{d,h} & \diff^2_{s\eta}f^{d,h} \\
\diff^2_{s\eta}f^{d,h} & \diff^2_{ss}f^{d,h}
\end{array}
\right)(\xi,0,u_1(\xi,0),u_2(\xi,0))>0,
\end{equation}
and
\begin{equation}\label{terzo}
\det \nabla^2_{\!\eta st}f^{d,h}(\xi,0,u_1(\xi,0),u_2(\xi,0))<0.
\end{equation}
By (\ref{primo}), (\ref{secondo}), and (\ref{terzo}),
we can conclude that the hessian matrix of
$f^{d,h}$ at $(\xi,0,u_1(\xi,0),u_2(\xi,0))$ is negative definite: 
both (\ref{disug}) and Claim~4 are proved.}

\bclaim{5}{For every $r>0$ and $h>0$,
there exists $\tilde{\varepsilon}>0$ with the property
that, if $\varepsilon \in (0, \tilde{\varepsilon})$, one can find 
$\delta\in (0, \varepsilon )$ so that 
\vspace{-.3cm}
\begin{eqnarray*}
& I(\xi, \eta, u_2(\xi,\eta),t )
\in B(r) + b_2^h(\xi, \eta, t), & \\
& I(\xi, \eta, s, u_1(\xi,\eta))
\in B(r) 
+ b_1^h(\xi, \eta, s), &
\end{eqnarray*} 
provided $V$ is small enough, for every $|t-u_2(\xi,\eta)|\leq \delta$,
$|s-u_1(\xi,\eta)|\leq \delta$.}
 
\eclaim{By the definition of $\phi$ in $A_6$, we obtain that 
$$I^{\xi}(\xi, \eta, u_2(\xi,\eta), t)
= 2(t-u_2(\xi,\eta))\diff_{\xi}u_2(\xi,\eta),$$ 
$$I^{\eta}(\xi, \eta, u_2(\xi,\eta), t) 
= 2(t-u_2(\xi,\eta))\diff_{\eta}u_2(\xi,\eta) 
- M (\varepsilon -M \eta)^{-1}(t-u_2(\xi,\eta))^2.$$ 
To get the claim, we need to prove that 
$$( 2(t-u_2)\diff_{\xi}u_2 )^2 
+ \left( -M (\varepsilon -M\eta )^{-1}(t-u_2)^2+ 
h(t-u_2)^2 + r \right)^2 
< r^2,$$ 
which is equivalent to  
$$( 2(t-u_2)\diff_{\xi}u_2 )^2 
+ \left(- M (\varepsilon -M\eta )^{-1} +h \right)^2(t-u_2)^4 
+ 2 r \left( -M (\varepsilon -M\eta )^{-1}+h \right) 
(t-u_2)^2 < 0.$$ 
The conclusion follows by remarking that, if $V$ is small enough,
the left-handside is less than 
$$\left(4(\diff_{\xi}u_2)^2 + 2hr-\frac{2Mr }{3\varepsilon}\right)\delta^2 
+o(\delta^2),$$ 
which is negative if $\varepsilon$ is sufficiently small.
The proof for $u_1$ is completely analogous.}

Let us conclude the proof of the step.
By Claim~4, we can find $h>0$ such that
(\ref{infpos}) is satisfied for $\varepsilon\in (0,
\overline{\varepsilon})$.
If we choose $r$ such that $2r < d/4$, by Claim~5 there exists
$\tilde{\varepsilon}>0$ such that
for every $\varepsilon\in (0, \tilde{\varepsilon})$ there is
$\delta\in (0,\varepsilon)$ so that
\begin{equation}\label{500}
\I{\eta}{s}{u_1(\xi,\eta)} + \I{\eta}{u_2(\xi,\eta)}{t} \in 
B(2r) + b_1^h(\xi,\eta, s) + b_2^h(\xi,\eta, t)
\end{equation}
for every $|s - u_1(\xi,\eta)|<\delta$, $|t - u_2(\xi,\eta)|<\delta$,
and $(\xi,\eta)\in V$.
If we take $\varepsilon \leq \min\{
\tilde{\varepsilon},  \overline{\varepsilon} \}$, then
by Claim~4 we have that the set
$$B(2r)+ (\emmee(\xi,\eta)\sin\thetat(\eta),\emmee
(\xi,\eta)\cos\thetat(\eta))
+ b_1^h(\xi,\eta, s) + b_2^h(\xi,\eta, t)$$
is contained in the ball centred at $(0,0)$ with radius $\gamma(\xi,\eta)$.
Some easy geometric considerations show that the relation
between $\thetae$ and $\thetat$ (see (\ref{angolo})) implies that
also the set 
\begin{equation}\label{insieme}
B(2r)+ (\emmee(\xi,\eta)\sin\thetae(\eta),\emmee
(\xi,\eta)\cos\thetae(\eta))
+ b_1^h(\xi,\eta, s) + b_2^h(\xi,\eta, t)
\end{equation}
is contained in the ball centred at $(0,0)$ with radius $\gamma(\xi,\eta)$, 
if the condition 
$$|b_1^h(\xi,\eta, s) + b_2^h(\xi,\eta, t)|< 2r$$
holds (to make this true, take 
$\delta$ and $V$ smaller if needed).
Since
$$\I{\eta}{s}{t}= 
\I{\eta}{s}{u_1(\xi,\eta)} + \I{\eta}{u_1(\xi,\eta)}{u_2(\xi,\eta)}
+ \I{\eta}{u_2(\xi,\eta)}{t},$$
by (\ref{500}), (\ref{pol1}), and (\ref{pol2}), it follows that
$\I{\eta}{s}{t}$ belongs to the set (\ref{insieme}), and then
to the ball centred at $(0,0)$
with radius $\gamma(\xi,\eta)$ for 
every $|s - u_1(\xi,\eta)|<\delta$, $|t - u_2(\xi,\eta)|<\delta$,
and $(\xi,\eta)\in V$. This concludes the proof of Step 1.

\bstep{2}{If $\varepsilon$ is sufficiently small and $\delta\in(0, \varepsilon)$,
condition (d) holds for
$|s-u_1(\xi,\eta)|\geq\delta$ or $|t-u_2(\xi,\eta)|\geq\delta$, and
$(\xi,\eta)\in V$, provided $V$ is small enough.}

\noindent Let us fix $\delta \in (0,\varepsilon)$ and set
$$m_1(\xi,\eta):= \max \{ |\I{\eta}{s}{t}|: \ 
u_1(\xi,\eta)- \varepsilon \leq s\leq t \leq u_2(\xi,\eta)+\varepsilon,
\,|t-u_2(\xi,\eta)|\geq \delta \}.$$
It is easy to see that the function $m_1$ is continuous. 
Let us prove that $m_1(\xi,0)< \gamma(\xi,0)=1$.

Fixed $(\xi,0)\in V$, $u_1(\xi,0)- \varepsilon \leq s\leq t \leq u_2(\xi,0)+\varepsilon$,
with $|t-u_2(\xi,0)|\geq \delta$, we can write
\begin{equation} \label{spezz1}
\I{0}{s}{t}  =  \I{0}{s}{u_1(\xi,0)} + \I{0}{u_1(\xi,0)}{u_2(\xi,0)} 
+ \I{0}{u_2(\xi,0)}{t}.
\end{equation}

\bclaim{6}{For every $r>0$ there exists $\varepsilon>0$ such that
$$\I{0}{u_2(\xi,0)}{t}\in B(r), \qquad \I{0}{s}{u_1(\xi,0)}\in B(r)$$
for $0<|s-u_1(\xi,0)|\leq \varepsilon$, $0<|t-u_2(\xi,0)|\leq \varepsilon$,
and $(\xi,0)\in V$.}

\eclaim{See the similar proof of Claim~5 above.}

By (\ref{spezz1}), (\ref{e/2}), (\ref{conde}), and Claim~6, it follows that
\begin{equation} \label{1000}
\I{0}{s}{t} \in (0,1) + \overline{B(r)} + B(r) = (0,1) + B(2r)
\end{equation}
for $0<|s-u_1(\xi,0)|\leq \varepsilon$, $\delta \leq |t-u_2(\xi,0)|\leq \varepsilon$.
If $r< 1/4$, the set $(0,1) + B(2r)$ is contained in the open ball centred
at $(0,0)$ with radius $1$.

It remains to study the case $|s-u_1|\geq \varepsilon$ and
the case $|t-u_2|\geq \varepsilon$. Let us consider the latter;
the former would be completely analogous.
We can write
\begin{eqnarray*}
\I{0}{s}{u_1(\xi,0)} & = & \I{0}{s\land (u_1(\xi,0) + \varepsilon)}{u_1(\xi,0)}
+ \I{0}{s\lor (u_1(\xi,0) + \varepsilon)}{u_1(\xi,0) + \varepsilon}, \\
\I{0}{u_2(\xi,0)}{t} & = & \I{0}{u_2(\xi,0)}{u_2(\xi,0) - \varepsilon}
+ \I{0}{u_2(\xi,0)- \varepsilon}{t}.
\end{eqnarray*}
Therefore, by (\ref{spezz1})
\begin{eqnarray}
\I{0}{s}{t} &  = & \I{0}{u_1(\xi,0)}{u_2(\xi,0)}  +
\I{0}{s\land (u_1(\xi,0) + \varepsilon)}{u_1(\xi,0)} + \nonumber \\
& & +
\I{0}{u_2(\xi,0)}{u_2(\xi,0) - \varepsilon} 
+ \I{0}{s\lor (u_1(\xi,0) + \varepsilon)}{t} -\nonumber\\ 
& & - \I{0}{u_1(\xi,0) + \varepsilon}{u_2(\xi,0)- \varepsilon}. \label{spezzamento}
\end{eqnarray}
If $-2\varepsilon(\diff_{\xi}u_1(\xi,0)+\diff_{\xi}u_2(\xi,0))\geq 0$, we define 
$$C := [0,-2\varepsilon(\diff_{\xi}u_1(\xi,0)+\diff_{\xi}u_2(\xi,0))]
{ \times }[0, 1-2\varepsilon M];$$
if $-2\varepsilon(\diff_{\xi}u_1(\xi,0)+\diff_{\xi}u_2(\xi,0))< 0$, we simply replace 
$[0,-2\varepsilon(\diff_{\xi}u_1(\xi,0)+\diff_{\xi}u_2(\xi,0))]$ by 
$[-2\varepsilon(\diff_{\xi}u_1(\xi,0)+\diff_{\xi}u_2(\xi,0)),0]$.
From the definition of $\phi$ in $A_3\cup A_4\cup A_5$, it follows that 
\begin{equation}\label{600}
\I{0}{u_1(\xi,0) + \varepsilon}{u_2(\xi,0)- \varepsilon}=
(-2\varepsilon(\diff_{\xi}u_1(\xi,0)+\diff_{\xi}u_2(\xi,0)),
1-2\varepsilon M)
\end{equation}
and 
\begin{equation}\label{601}
\I{0}{s}{t}\in C
\end{equation}
for $u_1(\xi,0)+\varepsilon\leq s\leq t\leq u_2(\xi,0)-\varepsilon$.
Let $D:=C-(-2\varepsilon(\diff_{\xi}u_1(\xi,0)+\diff_{\xi}u_2(\xi,0)),
1- 2\varepsilon M)$.
Since $I^{\eta}(\xi,0, u_2(\xi,0), u_2(\xi,0) - \varepsilon)= -M\varepsilon$,
from (\ref{spezzamento}), (\ref{e/2}), (\ref{conde}), Claim~6, (\ref{600}),
and (\ref{601}), we obtain 
\begin{eqnarray*}
\I{0}{s}{t} & \in & [(0,1)+\overline{B(r)}+B(r)]
\cap \{(x,y)\in \rdue : y< 1- \varepsilon M \} \, +D \\
& & = [(0,1)+B(2r)] \cap \{(x,y)\in \rdue : y< 1- \varepsilon M \}\, +D.
\end{eqnarray*}
If $r<1/4$ and if $\varepsilon$ is sufficiently small, the set 
$[(0,1)+B(2r)] \cap \{(x,y)\in \rdue : y< 1- \varepsilon M \}\, +D$ 
is contained in the open ball centred at $(0,0)$ with radius $1$ and this means that
$m_1(\xi,0)<\gamma(\xi,0)$. 

Analogously we define 
$$m_2(\xi,\eta):= \max \{ |\I{\eta}{s}{t}|: \ 
u_1(\xi,\eta)- \varepsilon \leq s\leq t \leq u_2(\xi,\eta)+\varepsilon,
\, |s-u_1(\xi,\eta)|\geq \delta \}.$$
Arguing as in the case of $m_1$, we can prove that $m_2$ is continuous and $m_2(\xi,0)<\gamma(\xi,0)$. By continuity, if $V$ is small enough, 
$m_1(\xi,\eta)<\gamma(\xi,\eta)$ and $m_2(\xi,\eta)<\gamma(\xi,\eta)$,
for every $(\xi,\eta)\in V$: Step 2 is proved.
\vspace{.5cm}

By Step 1 and Step 2, we conclude that, choosing $\varepsilon$ sufficiently small
and $n=n_{\varepsilon}$ (see (\ref{neps})),
condition (d) is true for 
$u_1(\xi,\eta)-\varepsilon\leq s, t\leq u_2(\xi,\eta)+\varepsilon$ 
and in fact for every $s,t\in\erre$, from the definition of $\phi$ in $A_1$ and $A_7$.
\qed

\

\section{The graph-minimality}

We start this section with a negative result: if the domain $\Omega$ is too large,
the Euler conditions do not guarantee
the graph-minimality introduced in Definition \ref{graphmin},
as the following counterexample (suggested by Gianni Dal Maso) shows.

\begin{prop}
Let $R$ be the rectangle $(1, 1+4l){ \times } (-l,l)$
and let
$$u(x,y):= 
\begin{cases}
x & \text{if $y\geq 0$,} \\
-x & \text{if $y<0$.}
\end{cases}$$
Then, $u$ satisfies the Euler conditions for the Mumford-Shah functional in
$R$, but it is not a local graph-minimizer in $R$ for $l$ large enough.
\end{prop}

\noindent {\sc Proof.}
The Euler conditions are obviously satisfied by $u$ in $R$.

Let $R_0$ be the rectangle $(0,4) { \times } (-1,0)$ and
let $w$ be any function in $H^1(R_0)$
such that $w(x,0)=x$ for $x\in (0,2)$, and
$w(x,y)=0$ for $(x,y)\in \partial R_0 \setminus ((0,4){\times} \{ 0 \})$.

The idea is to perturb $u$ by the rescaled function
$v(x,y):= l w(\frac{x-1}{l}, \frac{y}{l})$.
We define the perturbed function
$$\tilde{u}(x,y):=
\begin{cases}
x & \text{on $R_1\setminus T_{\varepsilon}$,} \\
-x + \eta\, (x-1) & \text{on $T_{\varepsilon}$,} \\
-x + \eta\, v(x,y) & \text{on $R_2$,}
\end{cases}$$
where $\eta$ is a positive parameter and the rectangles $R_1$, $R_2$, and the
triangle $T_{\varepsilon}$ are indicated in Fig.~1.
\begin{figure}[h!]
\begin{center} 
  \includegraphics[height=0.4\textheight]{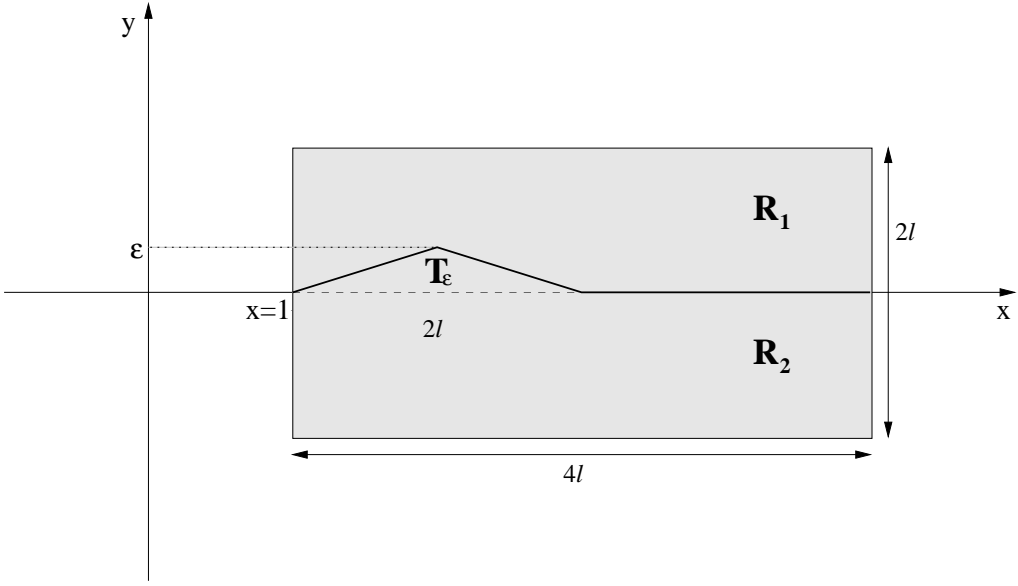}
  \caption{the regions $R_1$, $R_2$ and $T_{\varepsilon}$.}
\end{center}
\end{figure}
We want to show that, if we set $c:=\int_{R_0}|\nabla w(x,y)|^2 dx\,dy$,
for every $l>c$ and for every $\varepsilon_0$, $\eta_0 >0$
there exist $\varepsilon < \varepsilon_0$ and $\eta < \eta_0$
such that
$$\int_{R}|\nabla u(x,y)|^2 dx\,dy + {\cal H}^1(S_u) >
\int_{R}|\nabla \tilde{u}(x,y)|^2 dx\,dy + {\cal H}^1(S_{\tilde{u}}).$$

By definition, $\tilde{u}$ satisfies the boundary conditions.
Since by the construction of $v$ the function $\tilde{u}$ is continuous on the
interface between $T_{\varepsilon}$ and $R_2$, then
\begin{equation}\label{2003}
{\cal H}^1(S_u) - {\cal H}^1(S_{\tilde{u}})
= 2l - 2\sqrt{l^2 + \varepsilon^2} =
-\frac{\varepsilon^2}{l} + o(\varepsilon^2).
\end{equation}
On the triangle $T_{\varepsilon}$, we obtain
\begin{equation}\label{2001}
\int_{T_{\varepsilon}}|\nabla u(x,y)|^2 dx\,dy
- \int_{T_{\varepsilon}} |\nabla \tilde{u}(x,y)|^2dx\,dy
= 2l\varepsilon \eta - l\varepsilon \eta^2.
\end{equation}
Finally, since we have that 
$|\nabla \tilde{u}|^2 = 1 + \eta^2 |\nabla v|^2 
-2\eta\, \partial_x v$ in $R_2$, taking into account the boundary conditions of
$v$, we get
\begin{eqnarray}
\int_{R_2}|\nabla u(x,y)|^2 dx\,dy
- \int_{R_2} |\nabla \tilde{u}(x,y)|^2dx\,dy
& = & - \eta^2 \int_{R_2}|\nabla v(x,y)|^2 dx\,dy \nonumber \\
& = & - l^2 \eta^2 \int_{R_0}|\nabla w(x,y)|^2 dx\,dy. \label{2002}
\end{eqnarray}
In order to conclude,
by (\ref{2003}), (\ref{2001}), and (\ref{2002}), 
we have to show that for $l$ large we can choose $\varepsilon$
and $\eta$ arbitrarily close to $0$ such that
$$-\frac{\varepsilon^2}{l} - c l^2 \eta^2 +
2l\varepsilon \eta - l\varepsilon \eta^2 + o(\varepsilon^2) >0.$$
If we choose $\eta = \varepsilon/(cl)$,
then the equality above reduces to
$$-\frac{\varepsilon^2}{l} + \frac{\varepsilon^2}{c}
+ o(\varepsilon^2) >0,$$
which is true if $l>c$.
\qed

\subsection{Proof of Theorem \ref{teoglob}}

From the definition of $d$ and $N$ (see (\ref{defd}) and
Claim~3 in the proof of Theorem \ref{teoparz}) it follows that there
is an absolute constant $\tilde{c}>0$ (independent of $\Omega_0$, $\Omega$,
$\Gamma$, and $u$) such that 
\begin{equation}\label{defc}
\tilde{c}\,(1+ l^2(\Gamma)k^2(\Gamma)) >\frac{16}{d}.
\end{equation}
The absolute constant $c$, which appears in (\ref{sufficiente}), is defined
by 
\begin{equation}\label{assc}
c := \max \left\{ \tilde{c}, \frac{64}{\pi^2} \right\}.
\end{equation}

Actually, to avoid problems of boundary regularity, we shall work not exactly in $\Om$, 
but in a little bit larger set. Let $\Om'$ 
be a $\Gamma$-admissible set such that
$\Om \subset \subset \Om' \subset \subset \Om_0$, and 
$$\frac{\min_{i=1,2}  K(\Gamma\cap\Om' ,\Om_i')}
{1 + l^2(\Gamma\cap \Om') + l^2(\Gamma\cap \Om') k^2(\Gamma\cap
\Om')} > c  \sum_{i=1}^2 \|\partial_{\tau}
u_i\|^2_{C^1(\Gamma\cap\Om')},$$ 
where $\Om'_i$ denote the connected components of $\Om'\setminus \Gamma$.
This is possible by (\ref{sufficiente}) and
by the continuity properties of $K$.

The idea of the proof is to construct first a calibration $\varphi$
in a cylinder with base 
an open neighbourhood  of $\Gamma \cap \Om'$, and
then to extend $\varphi$ in a tubular neighbourhood
of $\graph{u}$. 

\begin{itemize}
\item {\it Construction of the calibration around $\Gamma$.}
\end{itemize}
We essentially recycle the construction of Theorem \ref{teoparz},
but we need to slightly modify the definition around
the graph of $u$,
in order to exploit condition (\ref{sufficiente}) and get the extendibility.

To define the calibration $\varphi(x,y,z)$ we use the same notation and the coordinate system 
$(\xi,\eta)$ on $U$ (open neighbourhood of $\Gamma\cap\Om'$) introduced in the
proof of Theorem \ref{teoparz}.  The vector field will be written as 
\begin{equation}\label{calibra22}
\varphi(x,y,z) = \frac{1}{\gamma^2 (\xi(x,y), 
\eta(x,y))}\phi(\xi(x,y), 
\eta(x,y), z),
\end{equation}
where $\phi$ can be represented by
$$\phi(\xi, \eta, z) = \phi^{\xi}(\xi, \eta, z) \ddxi + 
\phi^{\eta}(\xi, \eta, z) \ddeta
+ \phi^{z}(\xi, \eta, z) \ddz.$$

Given suitable parameters $\varepsilon >0$ and $\lambda >0$, 
we consider the following subsets of $V { \times } \erre$
\begin{eqnarray*}
A_1 &:= & \{(\xi,\eta,z)\in V{ \times } \erre : u_1(\xi,\eta) - 
	\varepsilon \, \vuu (\xi, \eta) < z < u_1(\xi,\eta) + 
	\varepsilon \, \vuu (\xi, \eta) \}, \\
A_2 & := & \{(\xi,\eta ,z)\in V{ \times } \erre : u_1(\xi,\eta) + \varepsilon \,
 	\vuu (\xi, \eta) < z <	u_1(\xi,\eta) + 2\varepsilon \}, \\ 
A_3 & := & \{(\xi,\eta,z)\in V{ \times } \erre : u_1(\xi,\eta) + 2 \varepsilon < z <
	\beta_1(\xi,\eta)  \}, \\
A_4 & := & \{(\xi,\eta,z)\in V{ \times } \erre : \beta_1(\xi,\eta) < z < 
	\beta_2(\xi,\eta) + 1/\lambda \}, \\
A_5 & := & \{(\xi,\eta,z)\in V{ \times } \erre : \beta_2(\xi,\eta)  + 1/\lambda
	< z < u_2(\xi,\eta) - 2 \varepsilon \}, \\
A_6 &:= & \{(\xi,\eta,z)\in V{ \times } \erre : u_2(\xi,\eta) - 2 \varepsilon 
	< z < u_2(\xi,\eta) - \varepsilon \, \vud (\xi, \eta ) \}, \\
A_7 & := & \{(\xi,\eta ,z)\in V{ \times } \erre : u_2(\xi,\eta) - \varepsilon \,
	\vud (\xi, \eta ) < z < u_2(\xi,\eta) + \varepsilon \, \vud (\xi,\eta )
\},   
\end{eqnarray*} 
where the functions $\vui$ are defined as
$$v_1(\xi, \eta) := 1 +M \eta, \;
v_2(\xi, \eta) := 1 -M \eta$$
with $M$ positive parameter such that
\begin{equation}\label{condicio}
c\, (1 + l^2(\Gamma\cap \Om') + l^2(\Gamma\cap \Om') k^2(\Gamma\cap
\Om'))  \sum_{j=1}^2 \|\partial_{\tau}
u_j\|^2_{C^1(\Gamma\cap\Om')} < M < \min_{j=1,2} K(\Gamma\cap\Om', \Omega_i'),
\end{equation}
while $\beta_1$ and $\beta_2$ are 
the solutions of the Cauchy problems
(\ref{betai}). Since we suppose $u_2>0$ on $V$, if
$\varepsilon$ is small enough,
while $\lambda$ is sufficiently large, then the sets $A_1, \ldots, A_7$
are nonempty and disjoint, provided $V$ is sufficiently small.
 
The vector 
$\phi(\xi,\eta,z)$
introduced in (\ref{calibra22}) will be written as
$$\phi(\xi,\eta,z)= (\phi^{\xi\eta}(\xi,\eta,z), \phi^z(\xi,\eta,z)),$$
where $\phi^{\xi\eta}$ is the two-dimensional vector
given by the pair $(\phi^{\xi}, \phi^{\eta})$.
We define $\phi(\xi,\eta,z)$
as follows:
%
%DEFINIZIONE DEL CAMPO
%
$$
\begin{cases}
\displaystyle \left( 2 \nablav  u_1 
- 2 \frac{u_1 - z}{\vuu }\nablav \vuu ,\,
\left|\nablav  u_1 -  \frac{u_1 - z}{\vuu}\nablav \vuu \right|^2 \right)
& \text{in $A_1$}, \\
\\
\displaystyle \left( 2 \nablav (u_1 
+ \varepsilon \vuu)   
- 2 \frac{u_1 +\varepsilon \vuu - 
z}{\vuup } \nablav \vuup , \, 
\left| \nablav ( u_1 
+ \varepsilon  \vuu)   
-  \frac{u_1 +\varepsilon \vuu - 
z}{\vuup} \nablav \vuup \right|^2
\right)
& \text{in $A_2$}, \\
\\
\displaystyle (0, \,\omega_1(\xi,\eta))  
& \text{in $A_3$}, \\
\\
\displaystyle  (\lambda \sigma(\xi,\eta) \nablav w  ,\, \mu)  & \text{in 
$A_4$}, \\
\\
\displaystyle(0, \, \omega_2(\xi,\eta))  
& \text{in $A_5$}, \\
\\
\displaystyle \left( 2\nablav ( u_2 
-  \varepsilon   \vud ) - 2\frac{u_2 -
\varepsilon \vud -
z}{\vudp} \nablav \vudp, \,
\left| \nablav ( u_2 
-  \varepsilon  \vud) - \frac{u_2 -
\varepsilon \vud -
z}{\vudp} \nablav \vudp \right|^2
\right)
& \text{in $A_6$}, \\
\\ 
\displaystyle \left( 2\nablav u_2 
- 2 \frac{u_2 -
z}{\vud}\nablav \vud, \,
\left| \nablav u_2 
-  \frac{u_2 -
z}{\vud}\nablav \vud \right|^2
\right)
& \text{in $A_7$},
\end{cases}$$
where $\nabla$ denotes the gradient with respect to the variables
$(\xi,\eta)$, the functions $\vuip$ are defined by
$$\tilde{v}_1(\xi, \eta):=2\varepsilon + M'\eta, \;
\tilde{v}_2(\xi, \eta):=2\varepsilon -M'\eta$$
while
$$\omega_i(\xi,\eta) := 
\varepsilon^2 \left( M +
M'\frac{\vui(\xi, \eta)}{\vuip(\xi, \eta)} \right)^2
- (\partial_{\xi}u_i(\xi,\eta))^2
- (\partial_{\eta}u_i(\xi,\eta))^2$$
for $i=1,2$, and for every $(\xi,\eta)\in V$; we take the constant $\mu$ 
sufficiently large in order to get the required inequality between 
the horizontal and the vertical components of
the field (see condition (b) of Section 2), and $M'$ so large that $\omega_i$
is positive in $V$, provided $V$ is small enough.
We define $w$ as the solution of the Cauchy problem
\begin{equation}\label{pbma}
\begin{cases}
\triangle w = 0, \\
\displaystyle w(\xi,0) = -\frac{4\varepsilon}{1-\varepsilon M'
-6\varepsilon^2 M }\int_0^{\xi}
n(s) (\diff_{\xi} u_1(s,0) + \diff_{\xi} u_2(s,0)) \,ds, \\
\diff_{\eta} w (\xi,0) = n(\xi),
\end{cases}
\end{equation}
where $n$ is a positive analytic function that
must be chosen in a suitable way.
We define
$$\sigma(\xi,\eta):= \frac{1}{n(q(\xi,\eta))} (1- \varepsilon M'
-6\varepsilon^2 M ),$$
where the function $q$ is constructed in the same way as in (\ref{qdef}).

Let us prove that for a suitable choice of the involved parameters  
the vector field is a calibration in a suitable neighbourhood $U$ of
$\Gamma\cap\Om'$, which is equivalent to prove that $\phi$ satisfies (a), (b),
(c), (d), and (e) of page~\pageref{calk}. The proof of conditions (a), (b),
(c), and (e) is the same of Theorem \ref{teoparz}. The proof  of (d) is
split again in two steps.

\bstep{1}{For a suitable choice of $\varepsilon$ and of the function $n$
(see (\ref{pbma}))
there exists $\delta >0$ such that condition (d) holds
for $|s-u_1(\xi,\eta)|<\delta$, $|t-u_2(\xi,\eta)|<\delta$, and
$(\xi,\eta)\in V$, provided $V$ is small enough.}

\noindent We essentially repeat the proof given in Theorem \ref{teoparz}:
Claims~1, 2, 3, and 4 are still valid with the same proof (up to the obvious
changes due to the different definition of $\phi$). Claim~5 must be modified
as follows.

\bclaim{5}{For $h=  \frac{64}{\pi^2}l^2(\Gamma)
\sum_{i=1}^{2}\| \diff_{\xi}u_i \|_{C^1(\Gamma\cap \Om')}^2$,
there exist $r\in(0, d/8)$ and $\tilde{\delta}>0$ such that for every
$\delta\in(0,\tilde{\delta})$
\begin{eqnarray*}
& I(\xi, \eta, u_2(\xi,\eta),t )
\in B(r) + b_2^h(\xi, \eta, t), & \\
& I(\xi, \eta, s, u_1(\xi,\eta))
\in B(r) 
+ b_1^h(\xi, \eta, s), &
\end{eqnarray*} 
provided $V$ is small enough, for every $|t-u_2(\xi,\eta)|\leq \delta$,
$|s-u_1(\xi,\eta)|\leq \delta$.}

\eclaim{Using the definition of $\phi$ in 
$A_7$, the claim is equivalent to prove 
$$( 2(t-u_2)\diff_{\xi}u_2 )^2 
+ \left(- M (1 -M \eta)^{-1} +h \right)^2(t-u_2)^4 
+ 2 r \left( -M (1 -M \eta)^{-1} +h \right) 
(t-u_2)^2 < 0;$$
note that for $a_1\in (0,1)$ the left-handside is less than
$$\left(4\sum_{i=1}^2\| \partial_{\xi} u_i \|_{C^1(\Gamma\cap\Om')}^2+
2hr-\frac{2r}{1+a_1}M \right)\delta^2+o(\delta^2),$$ 
provided $V$ is small enough.
To obtain the claim, it is sufficient to prove that
\begin{equation}\label{eqcla}
\frac{2}{r}\sum_{i=1}^2\| \partial_{\xi} u_i \|_{C^1(\Gamma\cap\Om')}^2
< \frac{1}{1+a_1}M - h.
\end{equation}
Since by (\ref{condicio}), (\ref{defc}), and (\ref{assc}) we can write
$$M =
\left(\frac{16+a_2}{d}+ \frac{64}{\pi^2}l^2(\Gamma\cap
\Om')\right)\sum_{i=1}^2\| \partial_{\xi} u_i \|_{C^1(\Gamma\cap\Om')}^2,$$
with $a_2>0$, the inequality (\ref{eqcla}) is equivalent to
$$\frac{2}{r} < \left( \frac{1}{1+a_1} -1 \right) \frac{64}{\pi^2}
l^2(\Gamma\cap \Om') +
\frac{16+a_2}{d}\frac{1}{1+a_1},$$
which is true if $a_1$ is sufficiently small and 
$r$ is sufficiently close to $d/8$.
The proof for $u_1$ is completely analogous.}

To conclude the proof of the step, let $r$ and $h$ be as in Claim~5. 
If we choose
$\varepsilon<\overline{\varepsilon}$ and $\delta \leq\min \{
\tilde{\delta}, \varepsilon \}$, by Claim~5 we have that 
\begin{equation}\label{800}
\I{\eta}{s}{u_1(\xi,\eta)} + \I{\eta}{u_2(\xi,\eta)}{t} \in 
B(2r) + b_1^h(\xi,\eta, s) + b_2^h(\xi,\eta, t)
\end{equation}
for every $|s - u_1(\xi,\eta)|<\delta$, $|t - u_2(\xi,\eta)|<\delta$,
and $(\xi,\eta)\in V$; since $h$ satisfies (\ref{noname}) and $2r<d/4$, we can
apply Claim~4 to deduce that the set 
$$B(2r)+
(\emmee(\xi,\eta)\sin\thetat(\eta),\emmee (\xi,\eta)\cos\thetat(\eta))
+ b_1^h(\xi,\eta, s) + b_2^h(\xi,\eta, t)$$
is contained in the ball centred at $(0,0)$ with radius $\gamma(\xi,\eta)$.
Some easy geometric considerations show that the relation
between $\thetae$ and $\thetat$ (see (\ref{angolo})) implies that
also the set 
\begin{equation}\label{insieme'}
B(2r)+ (\emmee(\xi,\eta)\sin\thetae(\eta),\emmee
(\xi,\eta)\cos\thetae(\eta))
+ b_1^h(\xi,\eta, s) + b_2^h(\xi,\eta, t)
\end{equation}
is contained in the ball centred at $(0,0)$ with radius $\gamma(\xi,\eta)$, 
if the condition 
$$|b_1^h(\xi,\eta, s) + b_2^h(\xi,\eta, t)|< 2r$$
holds (to make this true, take 
$\delta$ and $V$ smaller if needed).
Since
$$\I{\eta}{s}{t}= 
\I{\eta}{s}{u_1(\xi,\eta)} + \I{\eta}{u_1(\xi,\eta)}{u_2(\xi,\eta)}
+ \I{\eta}{u_2(\xi,\eta)}{t},$$
by (\ref{500}), it follows that
$\I{\eta}{s}{t}$ belongs to the set (\ref{insieme'}), and then
to the ball centred at $(0,0)$
with radius $\gamma(\xi,\eta)$ for 
every $|s - u_1(\xi,\eta)|<\delta$, $|t - u_2(\xi,\eta)|<\delta$,
and $(\xi,\eta)\in V$. This concludes the proof of Step 1.

\bstep{2}{If $\varepsilon$ is sufficiently small and $\delta\in(0, \varepsilon)$,
condition (d) holds for
$|s-u_1(\xi,\eta)|\geq\delta$ or $|t-u_2(\xi,\eta)|\geq\delta$, and
$(\xi,\eta)\in V$, provided $V$ is small enough.}

\noindent By using condition (\ref{condicio}), arguing 
as in the proof of Claim~5, we can prove the following claim.

\bclaim{6}{There exist $r<1/4$ and $\varepsilon>0$ such that
$$\I{0}{u_2(\xi,0)}{t}\in B(r), \qquad \I{0}{s}{u_1(\xi,0)}\in B(r)$$
for $0<|s-u_1(\xi,0)|\leq \varepsilon$, $0<|t-u_2(\xi,0)|\leq \varepsilon$,
and $(\xi,0)\in V$.}
\vspace{.1cm}

We can conclude the proof of Step 2 in the same way as in Theorem
\ref{teoparz}, with the minor changes due to the different definition
of the field.
\vspace{.5cm}

By Step 1 and Step 2, we conclude that, choosing $\varepsilon$ sufficiently small
and $n$ in a suitable way,
condition (d) is true for 
$u_1(\xi,\eta)-\varepsilon\leq s, t\leq u_2(\xi,\eta)+\varepsilon$. So,
$\varphi$ is a calibration.

\begin{itemize}
\item {\it Construction of the calibration around the graph of $u$.}
\end{itemize}
Now the matter is to extend the field in a tubular neighbourhood 
of the graph of $u$.
From now on, we reintroduce the Cartesian coordinates.

Let $\Gamma_i$ be the curve $\eta = (-1)^i k$, where
$k>0$. If $k$ is sufficiently small, 
for $i=1,2$ the curve $\Gamma_i$ 
connects two points of $\partial \Om'_i$,
divides $\Om'_i$ (and
then $\Om$) in two connected components, and
the normal vector $\nu_i$ to $\Gamma_i$
which points towards $\Gamma$ coincides with
$(-1)^{i+1} \nabla \eta / |\nabla \eta |$.
Set $U':=U\cap \{(x,y)\in \Om' : \, |\eta(x,y) |< k \}$
and $U'':= U' \cap \Om$.
Since $\| \nabla \eta \| =1$ on $\Gamma$, by (\ref{condicio}) we can suppose
that  
\begin{equation}\label{condicio'}
\frac{M}{1-M k} \max_{i=1,2} \|\nabla\eta\|_{L^{\infty}(\Gamma_i)}
< \min_{i=1,2} K(\Gamma_i,\Om'_i\setminus \overline{U'}).
\end{equation}
Chosen $\delta$ so small that $(\graph{u})_{\delta}\cap((U''\cap\Om_1)\times\erre)
\subset A_1$ and $(\graph{u})_{\delta}\cap((U''\cap\Om_2)\times\erre)\subset A_7$,
we define the vector field
$$\hat{\varphi}(x,y,z) = (\hat{\varphi}^{xy}(x,y,z),
\hat{\varphi}^z(x,y,z)) \in \rtre,$$
as follows:
$$\begin{cases}
\varphi(x,y,z) & \text{in $\{(x,y,z):(x,y)\in U'',\ u_1(x,y)-\delta<z<u_2(x,y)+\delta \}$,}\\
\\
\displaystyle \left( 2\nabla u - 2 \frac{u-z}{\vuut}\nabla \vuut,
\left| \nabla u -  \frac{u-z}{\vuut}\nabla \vuut \right|^2 \right) &
\text{in $(\graph{u})_{\delta}\cap (\Om_1\setminus U'')
{ \times } \erre$,} \\
\\
\displaystyle \left( 2\nabla u - 2 \frac{u-z}{\vudt}\nabla \vudt,
\left| \nabla u -  \frac{u-z}{\vudt}\nabla\vudt \right|^2 \right) &
\text{in $(\graph{u})_{\delta}\cap (\Om_2\setminus U'')
{ \times } \erre$.}	
\end{cases}$$
The function
$\vuit$ is the solution of the problem
\begin{equation}\label{estensione}
\min \left\{ \int_{\Om'_i \setminus \usec} |\nabla v|^2 dx\,dy
- \frac{M}{1 - M k} \int_{\Gamma_i}
|\nabla \eta|\, v^2 d{\cal H}^1:
\, v\in H^1(\Om'_i \setminus \usec), 
\, v|_{\partial(\Om'_i \setminus \usec)
\setminus \Gamma_i} = 1 \right\}.
\end{equation}
Let us show that the problem (\ref{estensione}) admits a solution.
If $\{ v_n \}$ is a minimizing sequence, then
\begin{equation}\label{limitata}
\sup_n \int_{\Omega_i'\setminus\overline{U}'}|\nabla v_n|^2 dx\, dy - 
\frac{M}{1-M k}\int_{\Gamma_i}|\nabla\eta|\, v_n^2\, d{\cal H}^1
< + \infty.
\end{equation}
We have only to show that $\{ v_n \}$ is bounded in $H^1(\Om'_i \setminus \usec)$.
If we put $\overline{v}_n:=v_n-1$, 
by (\ref{capacita}) for every $\tau\in (0,1)$
we have 
\begin{eqnarray}
\lefteqn{\int_{\Om_i'\setminus\usec}|\nabla v_n|^2 dx\, dy
= \int_{\Om_i'\setminus\usec}|\nabla \overline{v}_n|^2 dx\, dy}\nonumber\\
& = & \left( \int_{\Gamma_i}\overline{v}_n^2 d{\cal H}^1 \right)
\int_{\Om_i' \setminus \usec}\left|\nabla\left(\frac{
\overline{v}_n}{(\int_{\Gamma_i}\overline{v}_n^2  d{\cal
H}^1)^\frac{1}{2}}\right)\right|^2  dx\, dy\nonumber\\
& \geq & \left( \int_{\Gamma_i}(v_n -1)^2 d{\cal H}^1 \right)
K(\Gamma_i,\Om_i'\setminus\usec) \nonumber\\
& \geq & (1-\tau)K(\Gamma_i,\Om_i'\setminus\usec)\int_{\Gamma_i}v_n^2 d{\cal H}^1
+ K(\Gamma_i,\Om_i'\setminus\usec) \, \left(1
-\frac{1}{\tau}\right){\cal H}^1(\Gamma_i)\label{coercivo},
\end{eqnarray}
where we used Cauchy Inequality.
By (\ref{condicio'}), we can
choose $\tau$ so small that 
$$(1-\tau )K(\Gamma_i,\Om_i'\setminus\usec)>
\frac{M}{1-M k}\|\nabla\eta\|_{L^{\infty}(\Gamma_i)},$$
and substituting (\ref{coercivo}) in (\ref{limitata}), we obtain 
$$\sup_n\int_{\Gamma_i}v_n^2\, d{\cal H}^1<+\infty.$$
Using again (\ref{limitata}) and Poincar\'e Inequality, 
we conclude that $\{ v_n\}$ is actually bounded in $H^1(\Om'_i \setminus \usec)$.

The solution of (\ref{estensione}) satisfies 
\begin{equation}\label{vi}
\begin{cases}
\triangle\vuit =0 & \text{in $\Om'_i \setminus \usec$,} \\
\displaystyle \frac{\diff \vuit}{\diff \nu} = 
\frac{M}{1 - M k}|\nabla \eta| \vuit & \text{on $\Gamma_i$,} \\
\vuit=1 & \text{on $\partial(\Om'_i \setminus \usec)\setminus \Gamma_i$,}
\end{cases}
\end{equation}
and so, in particular, belongs to
$C^{\infty}(\overline{\Om_i \setminus U''})$.
By a truncation argument, it is easy to see that
$\vuit \geq 1$, so 
$\hat{\varphi}$ is well defined.

Since $\hat{\varphi}$ is a calibration
in $\{(x,y,z):(x,y)\in U'',\ u_1(x,y)-\delta<z<u_2(x,y)+\delta \}$, it remains to prove
only that the field is globally divergence free
in the sense of distributions
and that conditions (b), (c), (d)
are verified in  the regions
$(\graph{u})_{\delta}\cap (\Om_i\setminus U'')
{ \times } \erre$.
First of all, note that
by Lemma \ref{alberti} the field
$\hat{\varphi}$ is divergence free in the regions
$(\graph{u})_{\delta}\cap (\Om_i\setminus U'')
{ \times } \erre$, since it is constructed starting from the family
of harmonic functions $u(x,y) -t\vuit(x,y)$.
To complete the proof, we need to check that the 
normal components of the traces 
of $\varphi$ and of the extension field are equal on the surface of separation, i.e.,
\begin{equation}\label{eqnorm}
\varphi^{xy} \cdot \nu_i =
\left(2\nabla u - 2\frac{u-z}{\vuit} \nabla \vuit \right)\cdot
\nu_i\qquad\hbox{on }\Gamma_i, \end{equation}
where $\nu_i=(-1)^{i+1} \nabla \eta / |\nabla \eta |$.
Using the definition of $\varphi$, we obtain that
$$\varphi^{xy} \cdot \nu_i = \left( (-1)^{i+1} \diff_{\eta} u -
\frac{u-z}{1 -M k}M \right) |\nabla \eta|;$$
since $\nabla u \cdot \nu_i = (-1)^{i+1} \diff_{\eta} u|\nabla \eta|$,
the equality (\ref{eqnorm}) is equivalent to
$$\frac{M}{1 - M k} |\nabla \eta|=
\frac{1}{\vuit}\nabla \vuit \cdot \nu_i,$$
which is true by (\ref{vi}).

Conditions (b) and (c) are obviously satisfied,
while condition (d) is true if we take $\delta$ satisfying
$$\delta \leq \sup \left\{ \left(4|\nabla u| + 2 \frac{|\nabla
\vuit |}{\vuit }  \right)^{-1}: \,
(x,y)\in \Om_i \setminus U'', \, i=1,2 \right\}.$$
Therefore, with this choice of $\delta$, the vector field
$\hat{\varphi}$ is a calibration. 
\qed

\subsection{Some properties of $K(\Gamma,A)$}

In this subsection we investigate some qualitative properties 
of the quantity $K(\Gamma,A)$ 
and we shall compute it explicitly in a very particular case.
Let us start by a very simple result.

\begin{prop}\label{explo}
Let $\Gamma$ be a simple analytic curve and $\tilde{\Gamma}$
an extension of $\Gamma$, whose endpoints do not coincide with the endpoints
of $\Gamma$. If $\Gamma_{\delta}^{\pm}$ are the
two connected components of $\Gamma_{\delta}\setminus \tilde{\Gamma}$
(which are well defined if $\delta$ is sufficiently small), then
$$\lim_{\delta \to 0^+} K(\Gamma, \Gamma_{\delta}^{\pm})=+\infty.$$  
\end{prop}

\noindent {\sc Proof.} For convenience we set
$$W^{\pm}(\delta ):= \left\{ v\in H^1(\Gamma_{\delta}^{\pm}): \ \int_{\Gamma}
v^2 d\haus =1, \ v=0 \ \hbox{on} \ \partial (\Gamma_{\delta}^{\pm})
\setminus \Gamma \right\}.$$
Suppose by contradiction that there exists a sequence $\{ \delta_n
\}$ decreasing to $0$ such that
$\sup_n K(\Gamma, \Gamma_{\delta_n}^{+}) = c <+\infty$; this implies the
existence of a sequence $\{ v_n \}$ such that
$$v_n \in W^{+}(\delta_n) \qquad \hbox{and} \qquad 
\int_{\Gamma_{\delta_n}^{+}} |\nabla v_n (x,y) |^2 dx\, dy \leq c$$
for every integer $n$.
From now on, we regard $v_n$ as a function belonging to
$H^1(\Gamma_{\delta_1}^{+})$ which vanishes on $\Gamma_{\delta_1}^{+}
\setminus \Gamma_{\delta_n}^{+}$. By Poincar\'e Inequality it follows
immediately that $\{ v_n \}$ is bounded in $H^1(\Gamma_{\delta_1}^{+})$, and so
admits a weakly convergent subsequence $\{v_{n_k} \}$.
Let us call $v$ the limit of the subsequence; since for every $k$,
$v_{n_k}$ vanishes on $\Gamma_{\delta_1}^{+}
\setminus \Gamma_{\delta_{n_k}}^{+}$, then $v$ must vanish a.e.;
on the other hand, since $\int_{\Gamma}
v_{n_k}^2 d\haus =1$, by the compactness of the trace operator, we have that
$\int_{\Gamma} v^2 d\haus =1$, and this is clearly impossible.
\qed

\

We remark that by Theorem \ref{teoglob} and Proposition \ref{explo}, 
if $U_0$ is a neighbourhood of $\Gamma$ and $u\in SBV(U_0)$ satisfies the Euler
conditions in $U_0$ with $S_u=\Gamma$, then
there exists a neighbourhood $U$ of $\Gamma$ contained in $U_0$ such that
$u$ is a local graph-minimizer in $U$.
Actually, taking $U$ smaller if needed, by Theorem \ref{teoparz}
we get also the Dirichlet minimality.

\begin{prop}\label{char}
(Characterization of $K(\Gamma,A)$.)
Let $A$ be an open set with Lipschitz boundary and $\Gamma$ be
a subset of $\partial A$ with nonempty relative interior in $\partial A$.
The constant
$K(\Gamma,A)$ is the first eigenvalue of the problem
\begin{equation}\label{eigen}
\begin{cases}
\Delta u=0 & \text{on $A$,}\cr
\displaystyle \frac{\partial u}{\partial\nu}=\lambda u & \text{on $\Gamma$,}\cr
u=0 & \text{on $\partial A\setminus\Gamma$.}
\end{cases}
\end{equation}
Moreover, it is the unique eigenvalue with a positive eigenfunction.
\end{prop}

\noindent {\sc Proof.} 
If $u$ is a solution of (\ref{capacita}), 
then it is harmonic and there exists a Lagrange multiplier $\lambda$ such that
\begin{equation}\label{lagrange}
2\int_{A}\nabla u\cdot\nabla\varphi\, dx\, dy
=\lambda\int_{\Gamma}u\varphi\, d{\cal H}^1\qquad 
\forall \varphi\in C^{\infty}(A):\
\varphi=0\hbox{ on }\partial A\setminus\Gamma,
\end{equation}
which means, by Green Formula, that ${\frac{\partial u}{\partial\nu}}=\lambda u$
on $\Gamma$.
Using (\ref{lagrange}), one can easily see that $K(\Gamma,A)$ 
is in fact the minimal eigenvalue of (\ref{eigen})
and that it has a positive eigenfunction (indeed, if $u$ is a solution also
$|u|$ is). Let $u$ be a positive function belonging to the eigenspace of
$K(\Gamma,A)$ and $v$ another positive eigenfunction associated with the
eigenvalue $\mu$; by Green Formula we have
$$\int_{\Gamma} v\frac{\partial u}{\partial\nu}\, d{\cal H}^1
- \int_{\Gamma} u\frac{\partial v}{\partial\nu}\, d{\cal H}^1=0,$$
therefore
$$(K(\Gamma,A) - \mu)\int_{\Gamma} uv\, d{\cal H}^1=0.$$
Since both $u$ and $v$ are positive, from the last equality
it follows that $\mu = K(\Gamma,A)$. \hfill \qed

\

\begin{prop} 
If $A = (0,a) {\times } (0,b)$ and $\Gamma = (0,a)
{ \times } \{ 0 \}$, then 
\begin{equation}\label{rett}
K(\Gamma, A) = \displaystyle \frac{\pi}{a\tanh\left( \frac{\pi
b}{a}\right)}. 
\end{equation}
\end{prop}

\noindent {\sc Proof.} The function
$$v(x,y)= \sin\left(\frac{\pi}{a}x\right) \sinh\left(
\frac{\pi}{a}(b-y)\right)$$
is positive and satisfies (\ref{eigen}) with
$\lambda = \displaystyle \frac{\pi}{a\tanh\left( \frac{\pi b}{a}\right)}$.
Then, by Proposition \ref{char}, this quantity coincides
with $K(\Gamma, A)$. \qed

\

\begin{prop}\label{prev}
Let $g: [0, a_0] \to [0, +\infty)$ be a Lipschitz function and denote
the graph of $g$ by $\Gamma$. Given $0\leq a_1 < a_2 \leq a_0$ and $b>0$,
if we set $\Gamma(a_1,a_2):= \graph{g}|_{(a_1, a_2)}$ and
$$R(a_1,a_2,b):= \{ (x,y): x\in (a_1, a_2),\ y\in (g(x), g(x)+b) \},$$
then
$$\lim_{|a_2 - a_1| \to 0} K\left( \Gamma(a_1,a_2),R(a_1,a_2,b) \right)=
+\infty \qquad \hbox{uniformly with respect to $b$}.$$  
\end{prop}

\noindent {\sc Proof.} The idea is to transform the region $R(a_1,a_2,b)$
into the rectangle $(0,a_2 - a_1) { \times } (0,b)$ by a suitable
diffeomorphism in order to use (\ref{rett}).

Let $\psi : (0,a_2 - a_1) { \times } (0,b)\to R(a_1,a_2,b)$ be the map
defined by $\psi(x,y) = ( x + a_1, y + g(x + a_1) )$.
Let $v\in H^1(R(a_1,a_2,b))$ be such that $v=0$ on $\partial
R(a_1,a_2,b)\setminus\Gamma (a_1,a_2)$ and
\begin{equation}\label{forma}
\int_{\Gamma(a_1,a_2)} v^2 d{\cal
H}^1 = \int_{0}^{a_2 - a_1} v^2(\psi(x, 0))\sqrt{1+ (g'(x))^2}\, dx =1.
\end{equation}
If we call $\tilde{v}(x,y):= v(\psi(x,y))$, then $\tilde{v}\in H^1((0,a_2 -
a_1) { \times } (0,b))$, $\tilde{v}=0$ on the boundary of the rectangle
except $(0,a_2 - a_1) { \times } \{ 0 \}$, and by (\ref{forma}) there 
exists $\lambda > 0$ such that $\lambda^2 \leq \sqrt{1+ \|g'\|_{\infty}^2}$ and
$$\lambda^2 \int_{0}^{a_2 - a_1} \tilde{v}^2(x, 0)\, dx =1.$$
Therefore, since $J\psi\equiv 1$,
\begin{eqnarray*}
\int_{R(a_1,a_2,b)}|\nabla v(x,y)|^2dx\,dy  & = &
\int_{(0,a_2 - a_1) { \times } (0,b)} |\nabla v(\psi(x,y))|^2dx\,dy \\
& \geq & (1+\|g'\|_{\infty}+\|g'\|_{\infty}^2)^{-1}
\int_{(0,a_2 - a_1) { \times } (0,b)} |\nabla \tilde{v}(x,y)|^2dx\,dy \\
& \geq & \lambda^{-2}(1+\|g'\|_{\infty}+\|g'\|_{\infty}^2)^{-1} K
\left((0,a_2 - a_1) { \times } \{ 0 \}, (0,a_2 - a_1) { \times } (0,b)
\vphantom{\sum} \right)  \\ 
& \geq & (1+
\|g'\|_{\infty}^2)^{-3/2} \frac{\pi}{2(a_2 - a_1)\tanh\left( \frac{\pi b}{a_2 -
a_1}\right)},  
\end{eqnarray*}
where the last inequality follows by the estimate on $\lambda$ and
by (\ref{rett}). Since $v$ is arbitrary, using the fact that
$0 < \tanh t \leq 1$ for every $t>0$, we obtain that
$$K\left( \Gamma(a_1,a_2),R(a_1,a_2,b) \right) \geq 
(1+ \|g'\|_{\infty})^{-3/2}
\frac{\pi}{2(a_2 - a_1)};$$
so, the conclusion is clear. \qed

\

We have already remarked (see Proposition \ref{explo}) that
the graph-minimality is guaranteed in small neighbourhoods of the 
discontinuity set $\Gamma$. As consequence of Proposition \ref{prev},
we obtain that the graph-minimality holds also in the open sets, which are
narrow along the direction parallel to $\Gamma$ and may be very large along the
normal direction. This is made precise by the following corollary.

\begin{figure}[h]
\begin{center} 
  \includegraphics[height= 0.4\textheight, width= 0.4\textwidth]{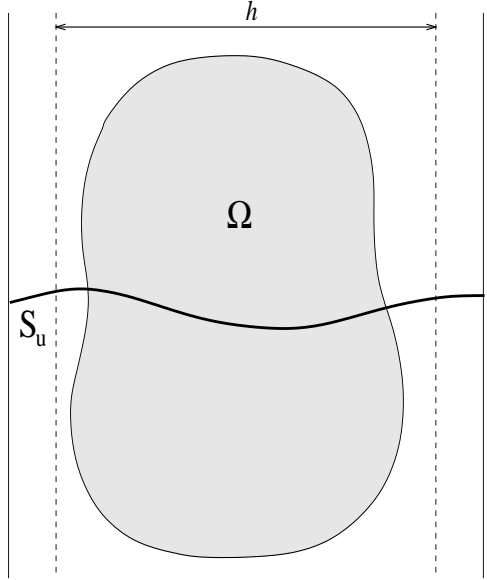}
\caption{if the thickness of $\Omega$ is less than $h$, then $u$ is a local
graph-minimizer in $\Omega$.}
\end{center}
\end{figure}

\begin{cor}\label{corolla}
Let $g$ be a positive function, analytic on $[0, a_0]$, that is
$g$ admits an analytic extension, and denote
the graph of $g$ by $\Gamma$.
For every $M>0$ there exists $h=h(M,\Gamma)$ 
such that, if $\Omega$ is $\Gamma$-admissible (see Definition \ref{gadm})
and $\Omega \subset (a_1, a_1 +h) { \times } \erre$
with $a_1\in [0, a_0 - h]$,
and if $u$ is a function in $SBV(\Omega)$ with $S_u=\Gamma\cap \Om$,
with different traces at every point of $\Gamma\cap \Om$, 
satisfying the Euler conditions in $\Omega$, and
$\sum_{i=1}^2 \| \diff_{\tau} u_i \|_{C^1(\Gamma\cap \Omega)}
\leq M$ (where $u_i$ is as above the restriction of $u$
to the connected component $\Omega_i$ of $\Omega \setminus
\Gamma$), then $u$ is a local graph-minimizer in
$\Omega$.  (see Fig.~2)
\end{cor}

\noindent {\sc Proof.} By Proposition \ref{prev} there exists
$h>0$ such that for every $a_1, a_2 \in [0, a_0]$ with
$0< a_2 - a_1 \leq h$ and for every $b>0$,
$$\frac{K(\Gamma(a_1, a_2), R(a_1, a_2, b))}
{1 + l^2(\Gamma ) + l^2(\Gamma)
k^2(\Gamma)} > c \,M^2.$$
If $\Omega \subset (a_1, a_1 +h) { \times } \erre$, then we can choose $b>0$
so large that, assuming that $\Omega_1$ is the upper component,
$\Omega_1 \subset R(a_1, a_1 + h, b)$.
Then by the monotonicity properties of $K(\Gamma,A)$, it follows that
$$\frac{K(\Gamma\cap \Omega , \Omega_1)}
{1 + l^2(\Gamma ) + l^2(\Gamma)
k^2(\Gamma)} > c \,M^2 \geq c \sum_{i=1}^2 \| \diff_{\tau} u_i
\|^2_{C^1(\Gamma\cap \Omega)}.$$
Applying the same argument to $\Omega_2$,
the conclusion follows from 
Theorem \ref{teoglob}.
\qed
%INSERT FIGURE 2

\

\section*{Acknowledgements}
We would like to thank Gianni Dal Maso for many helpful
discussions and for having suggested to us the study of this problem.

\end{document}